\documentclass[11pt,psamsfonts]{amsart}
\usepackage{amsmath}
\usepackage{amsthm}
\usepackage{amssymb}
\usepackage{amscd}
\usepackage{amsfonts}
\usepackage{amsbsy}
\usepackage{epsfig,afterpage}
\usepackage[dvips]{psfrag}
\usepackage[all]{xy}
\usepackage{pstcol}
\newcommand{\R}{\ensuremath{\mathbb{R}}}

\newcommand{\e}{\epsilon}

\newtheorem {theorem} {Theorem} 

\newtheorem {definition}  {Definition}
\newtheorem {remark}  {Remark}

\definecolor{red}{rgb}{1.,0.,0.}
\definecolor{blue}{rgb}{0.,0.,1.}
\definecolor{pink}{rgb}{1.,0.75,0.8}

\begin{document}

\title[Bifurcations of Non Smooth Vector Fields on $\R^2$] {Bifurcations
of Non Smooth Vector Fields on $\R^2$ by Geometric Singular
Perturbations}

\author[  T. de Carvalho, D.J. Tonon]
{ Tiago de Carvalho$^1$ and Durval Jos\'{e} Tonon$^2$}

\address{$^1$  IBILCE--UNESP, CEP 15054--000
S. J. Rio Preto, S\~ao Paulo, Brazil\\
$^2$ Universidade Federal de Goi\'{a}s, IME, CEP 74001-970 – Caixa
Postal 131, Goi\^{a}nia, GO, Brazil.}

\email{tiago@ibilce.unesp.br}

\email{djtonon@mat.ufg.br}


\date{}
\dedicatory{} \maketitle


\begin{abstract}
Our object of study is non smooth vector fields on $\R^2$. We apply
the techniques of geometric singular perturbations in non smooth
vector fields after regularization and a blow$-$up. In this way we
are able to bring out some results that bridge the space between
non$-$smooth dynamical systems presenting typical singularities and
singularly perturbed smooth systems.
\end{abstract}


\section{Introduction}

This work fits within the geometric study of Singular Perturbation
Pro\-blems expressed by vector fields on $\R^2$. We study the phase
portraits of certain non$-$smooth planar vector fields having a
curve $\Sigma$ as the dis\-con\-ti\-nuity set. We present some
results in the framework developed by Sotomayor and Teixeira in
\cite{ST} (and extended in \cite{LT-Regularization-1997}) and
establish a bridge between those systems and the fundamental role
played by the Geometric Singular Perturbation (abbreviated by GSP)
Theory.  This transition was introduced in papers like
\cite{Claudio-PR-Marco} and \cite{LST-Regularization-2006}, in
dimensions $2$ and $3$ respectively. Results in this context can be
found in \cite{LST}. We deal with non$-$smooth vector fields
pre\-sen\-ting structurally unstable con\-fi\-gu\-ra\-tions and we
prove that these structurally unstable con\-fi\-gu\-ra\-tions are
carried over the GSP Problem associated. Some good surveys about GSP
Theory are
\cite{DR} and \cite{F}, among others.\\

Let $\mathcal{U} \subseteq \R ^{2}$ be an open set and $\Sigma
\subseteq \mathcal{U}$ given by $\Sigma =f^{-1}(0),$ where
 $f:\mathcal{U} \longrightarrow \R$ is a smooth function having $0\in
\R$ as a regular value (i.e. $\nabla f(p)\neq 0$, for any $p\in
f^{-1}({0}))$. Clearly $\Sigma$ is the separating boundary of the
regions $\Sigma_+=\{q\in \mathcal{U} | f(q) \geq 0\}$ and
$\Sigma_-=\{q \in \mathcal{U} | f(q)\leq 0\}$. We can assume that
$\Sigma$ is represented, locally
around a point $q=(x,y)$, by the function $f(x,y)=x.$\\

Designate by $\mathfrak{X}^r$ the space of $C^r-$vector fields on a
compact set $K \subset \mathcal{U}$ endowed with the $C^r-$topology
with $r\geq 1$ large enough or $r=\infty$. Call
$\widetilde{\Omega}^r=\widetilde{\Omega}^r(K,f)$ the space of vector
fields $Z: K \setminus\Sigma \rightarrow \R ^{2}$ such that
$$
 Z(x,y)=\left\{\begin{array}{l} X(x,y),\quad $for$ \quad (x,y) \in
\Sigma_+,\\ Y(x,y),\quad $for$ \quad (x,y) \in \Sigma_-,
\end{array}\right.
$$
where $X=(f_1,g_1)$, $Y = (f_2,g_2)$ are in $\mathfrak{X}^r$.  The
trajectories of $Z$ are solutions of  ${\dot q}=Z(q),$
which has, in general, discontinuous right$-$hand side.\\

In what follows we will use the notation
\[X.f(p)=\left\langle \nabla f(p), X(p)\right\rangle \quad \mbox{ and }
 \quad Y.f(p)=\left\langle \nabla f(p), Y(p)\right\rangle. \]

We  distinguish the following regions on the discontinuity set
$\Sigma:$
\begin{itemize}
\item [$\blacktriangleright$]$\Sigma_1\subseteq\Sigma$ is the \textit{sewing region} if
$(X.f)(Y.f)>0$ on $\Sigma_1$ .
\item [$\blacktriangleright$]$\Sigma_2\subseteq\Sigma$ is the \textit{escaping region} if
$(X.f)>0$ and $(Y.f)<0$ on $\Sigma_2$.
\item [$\blacktriangleright$]$\Sigma_3\subseteq\Sigma$ is the \textit{sliding region} if
$(X.f)<0$ and $(Y.f)>0$ on $\Sigma_3$.
\end{itemize}

Consider $Z \in \widetilde{\Omega}^r.$ The \textit{sliding vector
field} associated to $Z$ is the vector field  $Z^s$ tangent to
$\Sigma_3$ and defined at $q\in \Sigma_3$ by $Z^s(q)=m-q$ with $m$
being the point of the segment joining $q+X(q)$ and $q+Y(q)$ such
that $m-q$ is tangent to $\Sigma_3$ (see Figure \ref{fig def
filipov}). It is clear that if $q\in \Sigma_3$ then $q\in \Sigma_2$
for $-Z$ and then we  can define the {\it escaping vector field} on
$\Sigma_2$ associated to $Z$ by $Z^e=-(-Z)^s$. The \textit{sewing
vector field}\index{vector field!sewing} associated to $Z$ is the
vector field $Z^{w}$ defined in $q \in \Sigma_{1}$ as an arbitrary
convex combination of $X(q)$ and $Y(q)$, i.e., $Z^{w}(q)=\lambda
X(q) + (1-\lambda) Y(q)$ where $\lambda \in [0,1]$. In what follows
we use the notation $Z^\Sigma$ for all these cases.

Let $\Omega^r=\Omega^r(K,f)$ be the space of vector fields $Z: K
\rightarrow \R^2$ such that $Z \in \widetilde{\Omega}^r$ and
$Z(q)=Z^{\Sigma}(q)$ for all $q \in \Sigma$. We write $Z=(X,Y),$
which we will accept to be multivalued in points of $\Sigma$. The
basic results of differential equations, in this context, were
stated by Filippov in
\cite{Fi}. Related theories can be found in \cite{K, ST, T}.\\

\begin{figure}[!h]
\begin{center}
\psfrag{A}{$q$} \psfrag{B}{$q + Y(q)$} \psfrag{C}{$q + X(q)$}
\psfrag{D}{} \psfrag{E}{\hspace{1cm}$Z^\Sigma(q)$}
\psfrag{F}{\hspace{.6cm}$\Sigma_2$} \psfrag{G}{} \epsfxsize=5.5cm
\epsfbox{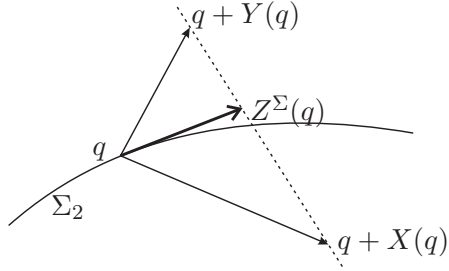} \caption{\small{Filippov's
convention.}} \label{fig def filipov}
\end{center}
\end{figure}

%
%
%
%

%

An approximation of the non$-$smooth vector field $Z = (X,Y)$ by a
$1-$parameter family $Z_{\epsilon}$ of smooth vector fields is
called an $\epsilon-$regularization of $Z$. We give the details
about this process in Section \ref{secao regularizacao}. A
transition function is used to average $X$ and $Y$ in order to get a
family of smooth vector fields that approximates $Z$. The main goal
of this process is to deduce certain dynamical properties of the
non$-$smooth dynamical system (abbreviated by NSDS) from the
regularized system. The regularization process developed by
Sotomayor and Teixeira produces a singular problem for which the
discontinuous set is a center manifold. Via a blow up we establish a
bridge between NSDS and the geometric singular perturbation theory.

Roughly speaking, the main results of this paper are the following:

\begin{theorem}\label{teorema boundary bifurcations}
Consider $Z(x,y)=Z_{\lambda}(x,y)=(X(x,y),Y_{\lambda}(x,y)) \in
\Omega^r$, where $\lambda\in (-1,1) \subset \R$, a non$-$smooth
planar vector field and $\Sigma$ identified with the $y-$axis. Let
the trajectories of $X$ be transverse to $\Sigma$ and $Y_{0}$
presenting either a hyperbolic saddle $q \in \Sigma$ or a hyperbolic
focus $q \in \Sigma$ or a $\Sigma-$cusp point $q$. Then there exists
a singular perturbation problem
\begin{equation}\label{eq pert sing}
\rho' = \alpha(r,\rho,y,\lambda) \, , \, y'=r
\beta(r,\rho,y,\lambda) \, ,
\end{equation}
with $r\geq 0$, $\rho \in (0,\pi)$, $y\in \Sigma$ and $\alpha$ and
$\beta$ of class $C^r$ such that the unfolding of \eqref{eq pert
sing} produces the same topological behaviors as the unfolding of
the corresponding normal forms of $Z_\lambda$ presented in
subsection \ref{subsecao formas normais}.
\end{theorem}

\begin{theorem}\label{teorema bif fold-fold}
Consider $Z(x,y)=Z_{\mu}(x,y)=(X_{\mu}(x,y),Y_{\mu}(x,y)) \in
\Omega^r$, where either $\mu=\lambda \in \R$ or
$\mu=(\lambda,\varepsilon) \in \R^2$, a non$-$smooth planar vector
field and $\Sigma$ identified with the $y-$axis. Consider that
$q=(x_q,y_q)\in \Sigma$ is a $\Sigma-$fold point of both $X_{\mu}$
and $Y_{\mu}$ when $\mu=0$ or $\mu=(0,0)$. Then there exists a
singular perturbation problem
\begin{equation}\label{eq pert sing 2}
\rho' = \alpha(r,\rho,y,\lambda) \, , \, y'=r
\beta(r,\rho,y,\lambda) \, ,
\end{equation}
with $r\geq 0$, $\rho \in (0,\pi)$, $y\in \Sigma$ and $\alpha$ and
$\beta$ of class $C^r$ such that the following statements holds:
\begin{description}
\item[(a)] For all small neighborhood $U$ of $q$ in $\Sigma$ the region
$(\Sigma_2 \cup \Sigma_3) \cap (U-\{ y_q \})$ is homeomorphic to the
slow manifold $\alpha(0,\rho,y,\lambda)=0$ of \eqref{eq pert sing 2}
where $y \in (U-\{ y_q \})$.

\item[(b)] The vector field $Z^{\Sigma}$, on $(\Sigma_2 \cup \Sigma_3) \cap (U-\{ y_q
\})$, and the reduced problem of \eqref{eq pert sing 2}, with $y \in
(U-\{ y_q \})$, are topologically equivalent.

\item[(c)] The slow manifold $\alpha(0,\rho,y,0)=0$ of \eqref{eq pert
sing 2}, where $y=y_q$, has just an horizontal component, i.e.,
$\alpha(0,\rho,y_q,0)=0$ can be identified with $\{ (\rho,y) \, | \,
\rho \in (0,\pi) \, , \, y=y_q \}$. Moreover, this configuration is
structurally unstable.

\end{description}

The unfolding of \eqref{eq pert sing 2} produces the same
topological behaviors as the unfolding of the corresponding normal
forms of $Z_\lambda$ presented in Table \ref{tabela} and in
\eqref{forma normal fold-fold eliptica}.
\end{theorem}

Observe that Theorem \ref{teorema bif fold-fold} generalize the
Theorem $1.1$ of \cite{LST}, because here we allow that
$X.f(q)=Y.f(q)=0$.

The paper is organized as follows: in Section \ref{preliminares} we
give the basic theory about Non$-$Smooth Vector Fields on the Plane,
in Section \ref{secao regularizacao} we give the theory about the
regularization process,  in Section \ref{secao perturbacao singular}
we present the GSP Theory, in Section \ref{secao boundary
bifurcations} we present the singularities treated in Theorem
\ref{teorema boundary bifurcations} and give its normal forms,  in
Section \ref{secao fold fold} we present the singularities treated
in Theorem \ref{teorema bif fold-fold} and give its normal forms and
in Section \ref{secao conclusao} we prove Theorems \ref{teorema
boundary bifurcations} and \ref{teorema bif fold-fold}.

\section{Preliminaries}\label{preliminares}

We say that $q\in\Sigma$ is a \textit{$\Sigma-$regular point} if
\begin{itemize}
\item [(i)] $X.f(q)Y.f(q)>0$
or
\item [(ii)] $X.f(q)Y.f(q)<0$ and $Z^{\Sigma}(q)\neq0$ (that is $q\in\Sigma_2\cup\Sigma_3$ and it is not a singular
point of $Z^{\Sigma}$).\end{itemize}

The points of $\Sigma$ which are not $\Sigma-$regular are called
\textit{$\Sigma-$singular}. We distinguish two subsets in the set of
$\Sigma-$singular points: $\Sigma^t$ and $\Sigma^p$. Any $q \in
\Sigma^p$ is called a \textit{pseudo equilibrium of $Z$} and it is
characterized by $Z^{\Sigma}(q)=0$. Any $q \in \Sigma^t$ is called a
\textit{tangential singularity} and is characterized by
$Z^{\Sigma}(q) \neq 0$ and
$X.f(q)Y.f(q) =0$ ($q$ is a contact point of $Z^{\Sigma}$).\\

A tangential singularity $q\in\Sigma^t$ is a \textit{$\Sigma-$fold
point} of $X$ if $X.f(q)=0$ but $X^{2}.f(q) = X.(X.f)(q)\neq0.$
Moreover, $q\in\Sigma$ is a \textit{visible} (resp. {\it invisible})
\textit{$\Sigma-$fold point} of $X$
 if $X.f(q)=0$ and $X^{2}.f(q)> 0$
(resp. $X^{2}.f(q)< 0$). We say that  $q\in\Sigma^t$ is a
\textit{$\Sigma-$cusp point} of $X$  if $X.f(q)=0$, $X^{2}f(q)=0$
but $X^{3}f(q)\neq0$. Moreover, $q\in\Sigma$ is a \textit{natural}
(resp. {\it inverse}) \textit{$\Sigma-$cusp point}  of $X$ if
$X.f(q)=0$, $X^{2}.f(q)=0$ and $X^{3}.f(q)> 0$ (resp. $X^{3}.f(q)<
0$).\\ 

A pseudo equilibrium $q \in \Sigma^p$ is a \textit{$\Sigma-$saddle}
provided one of the following condition is satisfied: (i)
$q\in\Sigma_2$ and $q$ is an attractor for $Z^{\Sigma}$ or (ii)
$q\in\Sigma_3$ and $q$ is a repeller for $Z^{\Sigma}$. A pseudo
equilibrium $q\in\Sigma^p$ is a $\Sigma-$\textit{repeller} (resp.
$\Sigma-$\textit{attractor}) provided $q\in\Sigma_2$ (resp. $q \in
\Sigma_3$) and $q$ is a repeller (resp. attractor) equilibrium point
for $Z^{\Sigma}$.

\section{Regularization}\label{secao regularizacao}

In this section we present the concept of $\epsilon-$regularization
of non$-$smooth vector fields. It was introduced by Sotomayor and
Teixeira in \cite{ST}. The regularization  gives the mathematical
tool to study the stability of these systems, according with the
program introduced by Peixoto. The method consists in the analysis
of the regularized vector field which is  a smooth approximation of
the non$-$smooth vector field. Using this process we get a
$1-$parameter family of vector fields $Z_{\e}\in \mathfrak{X}^r$
such that for each $\e_0 > 0$ fixed it satisfies that:
\begin{itemize}
\item [(i)]   $Z_{\e_0}$ is equal to $X$ in all
points of $ \Sigma_+$ whose  distance to $\Sigma$ is bigger than
$\e_0;$
\item [(ii)] $Z_{\e_0}$ is equal to $Y$ in all points of $ \Sigma_-$ whose distance to $\Sigma$ is bigger than $\e_0$.
\end{itemize}

\begin{definition}
 A $C^\infty$ function $\varphi:\R \longrightarrow \R$
is a transition function if $\varphi(x)=-1$ for $x\leqslant -1$,
$\varphi(x)=1$ for $x\geqslant 1$ and $\varphi'(x)>0$ if
$x\in(-1,1).$ The $\epsilon-$regularization of $Z=(X,Y)$ is the
$1-$parameter family $Z_{\e}\in \mathfrak{X}^r$ given by
\[
Z_{\e}(q)=\left(
\dfrac{1}{2}+\dfrac{\varphi_{\e}(f(q))}{2}\right)X(q) +\left(
\dfrac{1}{2}-\dfrac{\varphi_{\e}(f(q))}{2}\right) Y(q).
\]
with $\varphi_{\e}(x)=\varphi(x/\e),$ for $\e>0.$
\end{definition}


\section{Singular Perturbations}\label{secao perturbacao singular}



\begin{definition} Let $U\subseteq \R^2$ be an open subset and take $\e\geqslant 0$.
A singular perturbation problem in $U$ (SP$-$Problem) is a
differential system which can be  written like
\begin{equation}
\label{fast} x'=dx/d\tau=l(x,y,\e),\quad  y'=dy/d\tau=\e m(x,y,\e)
\end{equation} or equivalently, after the time re$-$scaling $t=\e\tau$
\begin{equation}
\label{slow} \e{\dot x}=\e dx/dt=l(x,y,\e),\quad {\dot y}=dy/dt=
m(x,y,\e),
\end{equation}
with $(x,y)\in U$ and $l,m$ smooth in all variables.
\end{definition}

The understanding of the  phase portrait of the vector field
associated to a SP$-$problem is the main goal of the
\textit{geometric singular perturbation theory} (GSP$-$theory). The
techniques of GSP$-$theory can be used to obtain information on the
dynamics of (\ref{fast}) for small values of $\e>0,$ mainly in
searching limit cycles. System (\ref{fast}) is called the
\textit{fast system}, and (\ref{slow}) the \textit{slow system} of
SP$-$problem. Observe that for $\e >0$ the phase portraits of the
fast and the slow systems coincide. For $\epsilon =0,$ let
$\mathcal{S}$ be the set
\[
 \label{SM} \mathcal{S}=\left\lbrace (x,y):f(x,y,0)=0\right\rbrace
\] of all singular points of (\ref{fast}). We call
$\mathcal{S}$ the slow manifold of the singular perturbation problem
and it is important to notice that equation (\ref{slow}) defines  a
dynamical system, on $\mathcal{S}$, called the {\it reduced
problem}:
\[
\label{reduced} f(x,y,0)=0 ,\quad {\dot y}= g(x,y,0).
\]
Combining results on the dynamics of these two limiting problems,
with $\e = 0$, one obtains information on the dynamics of $X_\e$ for
small values of $\e$. We refer to \cite{F} for an introduction to
the general theory of singular perturbations. Related problems can
be seen in \cite{BST}, \cite{DR} and \cite{S}.\\ 


\section{Boundary  Bifurcations}\label{secao boundary bifurcations}

Consider $Z=(X,Y)$. In this section we assume that the trajectories
of the smooth vector field $X$ is transversal to $\Sigma$ and that
$Y$ has either a hyperbolic saddle or a hyperbolic focus or
$\Sigma-$cusp point in $\Sigma$. This con\-fi\-gu\-ra\-tion is
clearly structurally unstable. We present here its normal forms and
unfoldings.

\subsection{Codimension One Normal Forms}\label{subsecao formas
normais}

Take $\Sigma$ as the $y-$axis, i.e., \linebreak$f(x,y)=x$ and
consider the parameter $\lambda \in (-1,1)$.

\begin{itemize}

\item \textit{Regular$-$saddle}: Assume that $X$ is transversal to $\Sigma$ and that $Y$ has  a
hyperbolic saddle in $\Sigma$. The following normal form generically
unfolds this configuration.
$$
 Z(x,y)=Z_{\lambda}(x,y)=\left\{\begin{array}{ll}
X(x,y) = \left(
           \begin{array}{c}
             1 \\
             1 \\
           \end{array}
         \right)
,& $for$ \quad (x,y) \in \Sigma_+, \\ Y_{\lambda}(x,y)= \left(
           \begin{array}{c}
             -y \\
             -x - \lambda\\
           \end{array}
         \right),& $for$ \quad
(x,y) \in \Sigma_-.
\end{array}\right.
$$

\item \textit{Regular$-$focus}: Assume that $X$ is transversal to $\Sigma$ and that $Y$ has  a
hyperbolic focus in $\Sigma$. The following normal form generically
unfolds this configuration.

$$
 Z(x,y)=Z_{\lambda}(x,y)=\left\{\begin{array}{ll}
X(x,y) = \left(
           \begin{array}{c}
             1 \\
             1 \\
           \end{array}
         \right)
,&$for$ \quad (x,y) \in \Sigma_+,\\ Y_{\lambda}(x,y)= \left(
           \begin{array}{c}
             x + y+ \lambda \\
             -x + y- \lambda \\
           \end{array}
         \right),& $for$ \quad
(x,y) \in \Sigma_-.
\end{array}\right.
$$

\item \textit{Regular$-$cusp}: Assume that $X$ is transversal to $\Sigma$ and that $Y$ has  a
$\Sigma-$cusp point. The following normal form generically unfolds
this configuration.

$$
 Z(x,y)=Z_{\lambda}(x,y)=\left\{\begin{array}{ll}
X(x,y) = \left(
           \begin{array}{c}
             1 \\
             1 \\
           \end{array}
         \right)
,& $for$ \quad (x,y) \in \Sigma_+,\\ Y_{\lambda}(x,y)= \left(
           \begin{array}{c}
             -y^2+ \lambda \\
             1 \\
           \end{array}
         \right),& $for$ \quad
(x,y) \in \Sigma_-,
\end{array}\right.
$$

\end{itemize}

\subsection{Regular$-$saddle Bifurcation}

Consider the regular$-$saddle normal form given in the previous
subsection.
The regularized vector field becomes
$$\begin{array}{lcl}
   \dot{x} & = & \dfrac{1-y}{2} +
\varphi \left( \frac{x}{\epsilon} \right)\dfrac{1+y}{2} \, ,
\\
   \dot{y}  & =  & \dfrac{1-x-\lambda}{2} + \varphi \left(
\frac{x}{\epsilon} \right)\dfrac{1+x+\lambda}{2}.\end{array}
$$
where $\varphi(x / \epsilon)$ is the transition function. Making the
polar blow up

\begin{equation} x=r \cos \theta \hspace{1cm}\mbox{ and }\hspace{1cm} \epsilon = r \sin
\theta, \label{blow-up polar}\end{equation} where $\theta = (2 \rho
+ \pi)/4,$ we obtain

\[
\begin{array}{lcl}
  r \dot{\theta}  & = & -\sin \theta \left(   \dfrac{1-y}{2} +
\varphi ( \cot \theta )\dfrac{1+y}{2} \right) \, , \\
  \dot{y} &= & \dfrac{1-r \cos \theta-\lambda}{2} + \varphi (
\cot \theta)\dfrac{1+r \cos \theta+\lambda}{2}.
\end{array}
\]

\begin{remark}\label{remark theta rho}
We use the new variable $\theta = (2 \rho + \pi)/4$, with  $\rho \in
(0, \pi)$, because in the coordinates $(r,\theta)$ the transition
function $\varphi$ is constant when $\theta \in (0,\pi/4) \cup (3
\pi/4,\pi)$. So, in the next Figures in the text, we can draw the
slow manifold and slow dynamics with $\rho \in (0,\pi)$.
\end{remark}

In the blowing up locus $r=0$ the fast dynamics is determined by the
system
 \[ \theta'   = -\sin \theta \left(   \dfrac{1-y}{2} +
\varphi ( \cot \theta )\dfrac{1+y}{2} \right) \, , \quad y'=0 \, ;
\] and the slow dynamics on the slow manifold is determined by the
reduced system
\[    \dfrac{-1+y}{2} +
\varphi ( \cot \theta )\dfrac{-1-y}{2}   = 0 \, , \quad \dot{y}=
\dfrac{1-\lambda}{2} + \varphi ( \cot \theta)\dfrac{1+\lambda}{2}.
\]

We remark that the slow manifold is implicitly defined by $(-1+y)/2
+ \varphi (\cot \theta) (-1-y)/2 = 0$ and do not depends on the
parameter $\lambda$ (see Figure \ref{fig pert singular sem nada e
invisivel}). Moreover, $y(\theta)$ defined in this way is such that
$$ \displaystyle\lim_{\theta \longrightarrow \frac{\pi}{4}}y(\theta)
= + \infty \quad \mbox{and} \quad \displaystyle\lim_{\theta
\longrightarrow \frac{3\pi}{4}}y(\theta) = 0. $$

\begin{figure}[!h]
\begin{minipage}{0.485\linewidth}
\psfrag{A}{$X$} \psfrag{B}{$Y$} \psfrag{C}{$y$} \psfrag{D}{$0$}
\psfrag{E}{$\frac{\pi}{4}$} \psfrag{F}{$\frac{\pi}{2}$}
\psfrag{G}{$\frac{3 \pi}{4}$}\psfrag{H}{$\pi$}\psfrag{I}{$x$}
\epsfxsize=6cm \epsfbox{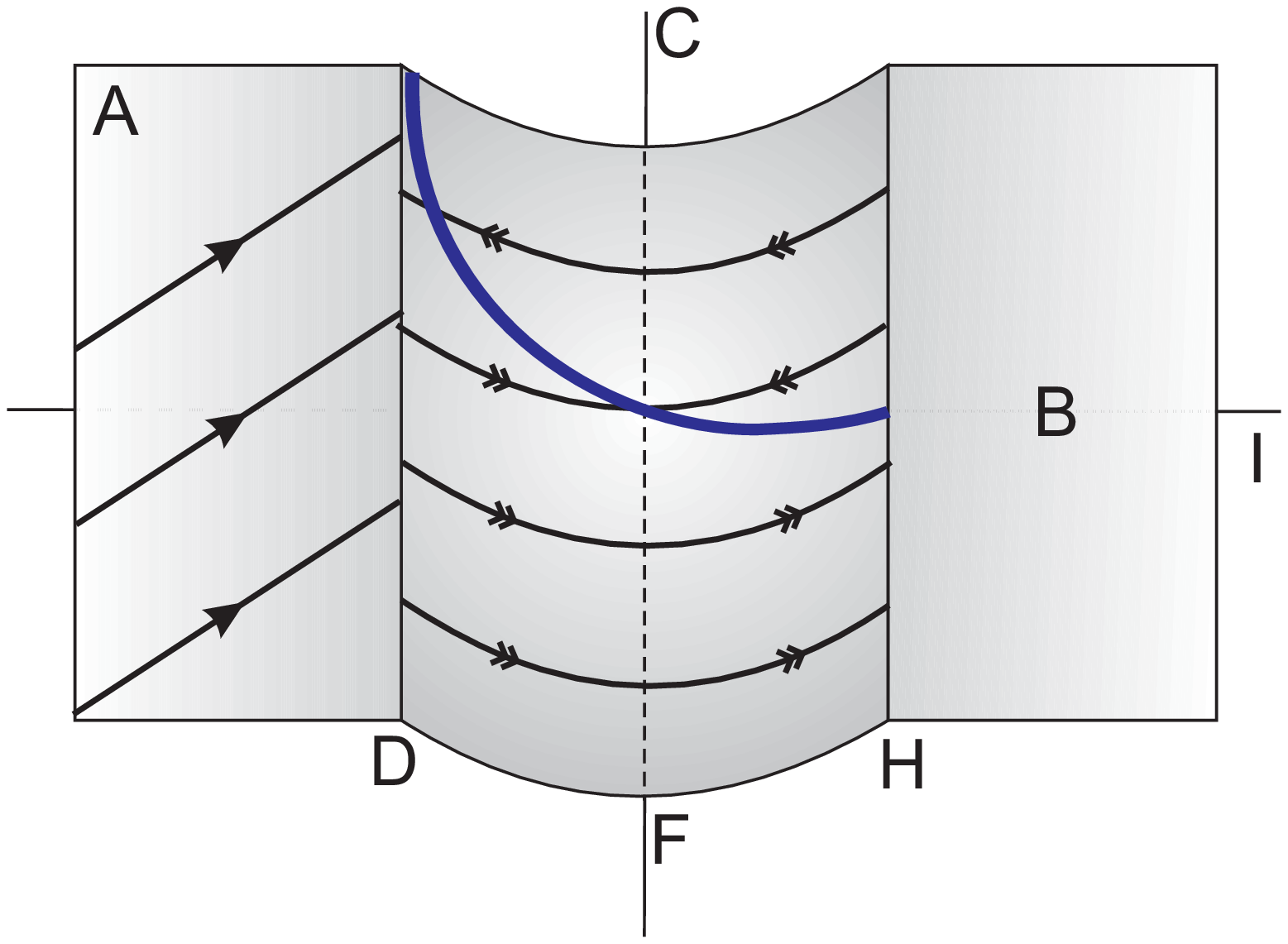}
\end{minipage} \hfill
\begin{minipage}{0.485\linewidth}
\psfrag{A}{$X$} \psfrag{B}{$Y$} \psfrag{C}{$y$} \psfrag{D}{$0$}
\psfrag{E}{$\frac{\pi}{4}$} \psfrag{F}{$\frac{\pi}{2}$}
\psfrag{G}{$\frac{3 \pi}{4}$}\psfrag{H}{$\pi$}\psfrag{I}{$x$}
\epsfxsize=6cm \epsfbox{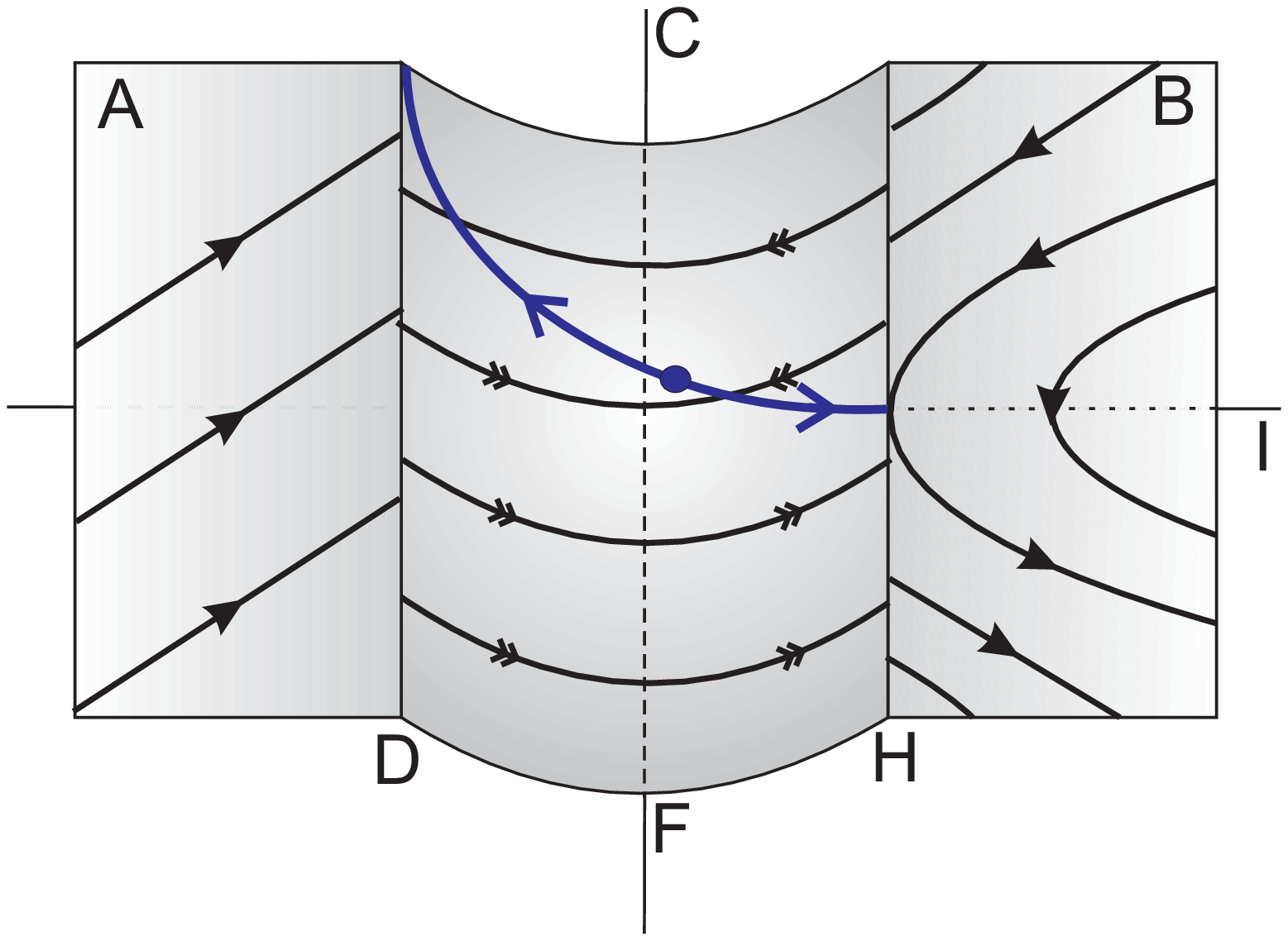}
\end{minipage}
\caption{\small{In this figure is pictured the slow manifold to the
left and the case $\lambda<0$ at the right. In both we consider
$\rho \in (0,\pi)$.}}\label{fig pert singular sem nada e invisivel}
\end{figure}

By other hand, the dynamics on the slow manifold depends on
$\lambda$. In fact, if either $\lambda>0$ or if $\lambda=0$ then
$\dot{y}\neq 0$ (see Figure \ref{fig pert singular sela sobre a
descontinuidade e visivel}) and if $\lambda>0$ so $\dot{y}$ has an
unique repeller critical point $P$   (see Figure \ref{fig pert
singular sem nada e invisivel}) given implicitly by the equation
$\varphi(\cot \theta)= (-1+\lambda)/(1 + \lambda)$.

\begin{figure}[!h]
\begin{minipage}{0.485\linewidth}
\psfrag{A}{$X$} \psfrag{B}{$Y$} \psfrag{C}{$y$} \psfrag{D}{$0$}
\psfrag{E}{$\frac{\pi}{4}$} \psfrag{F}{$\frac{\pi}{2}$}
\psfrag{G}{$\frac{3 \pi}{4}$}\psfrag{H}{$\pi$}\psfrag{I}{$x$}
\epsfxsize=6cm \epsfbox{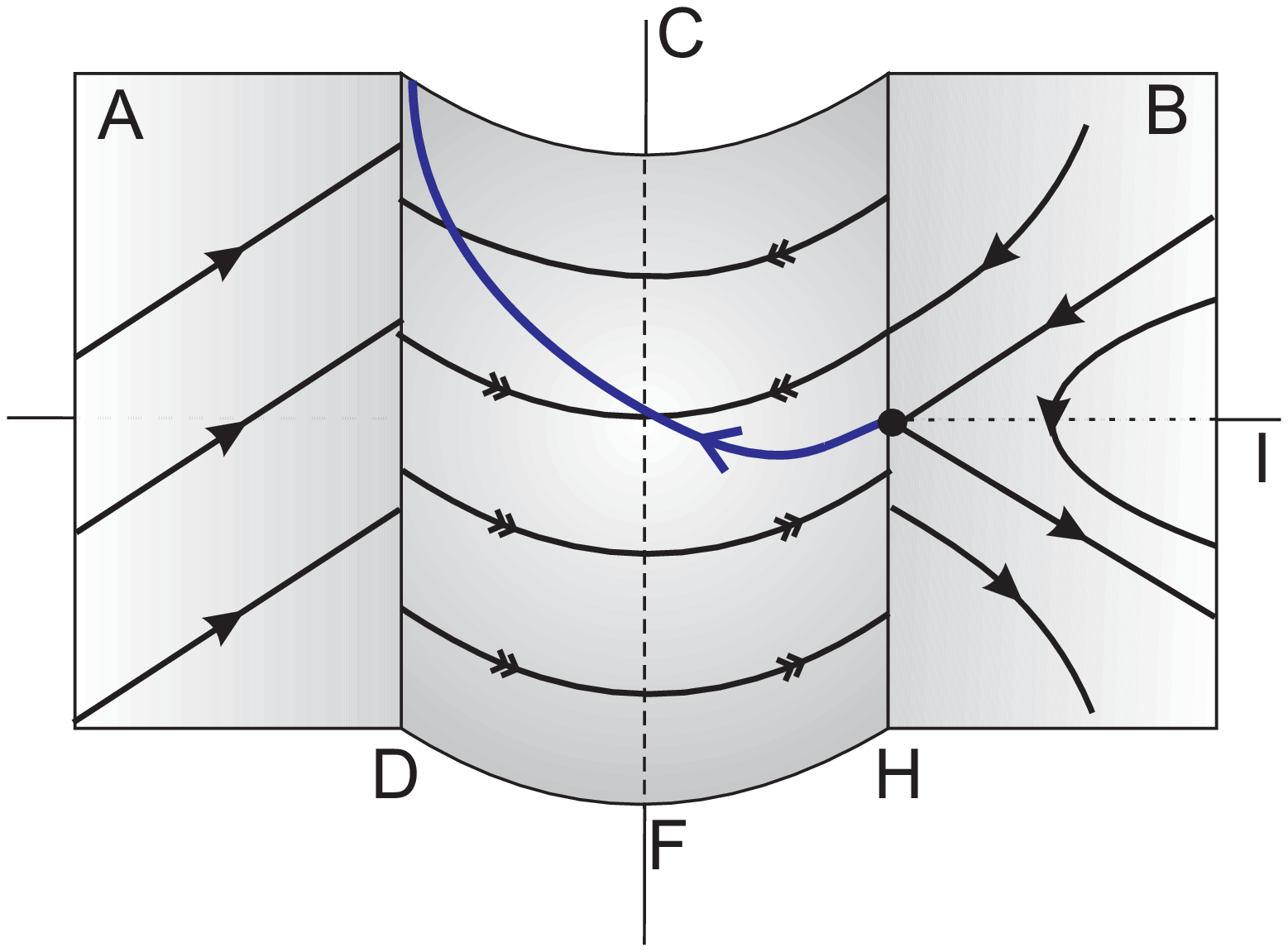}
\end{minipage} \hfill
\begin{minipage}{0.485\linewidth}
\psfrag{A}{$X$} \psfrag{B}{$Y$} \psfrag{C}{$y$} \psfrag{D}{$0$}
\psfrag{E}{$\frac{\pi}{4}$} \psfrag{F}{$\frac{\pi}{2}$}
\psfrag{G}{$\frac{3 \pi}{4}$}\psfrag{H}{$\pi$}\psfrag{I}{$x$}
\epsfxsize=6cm \epsfbox{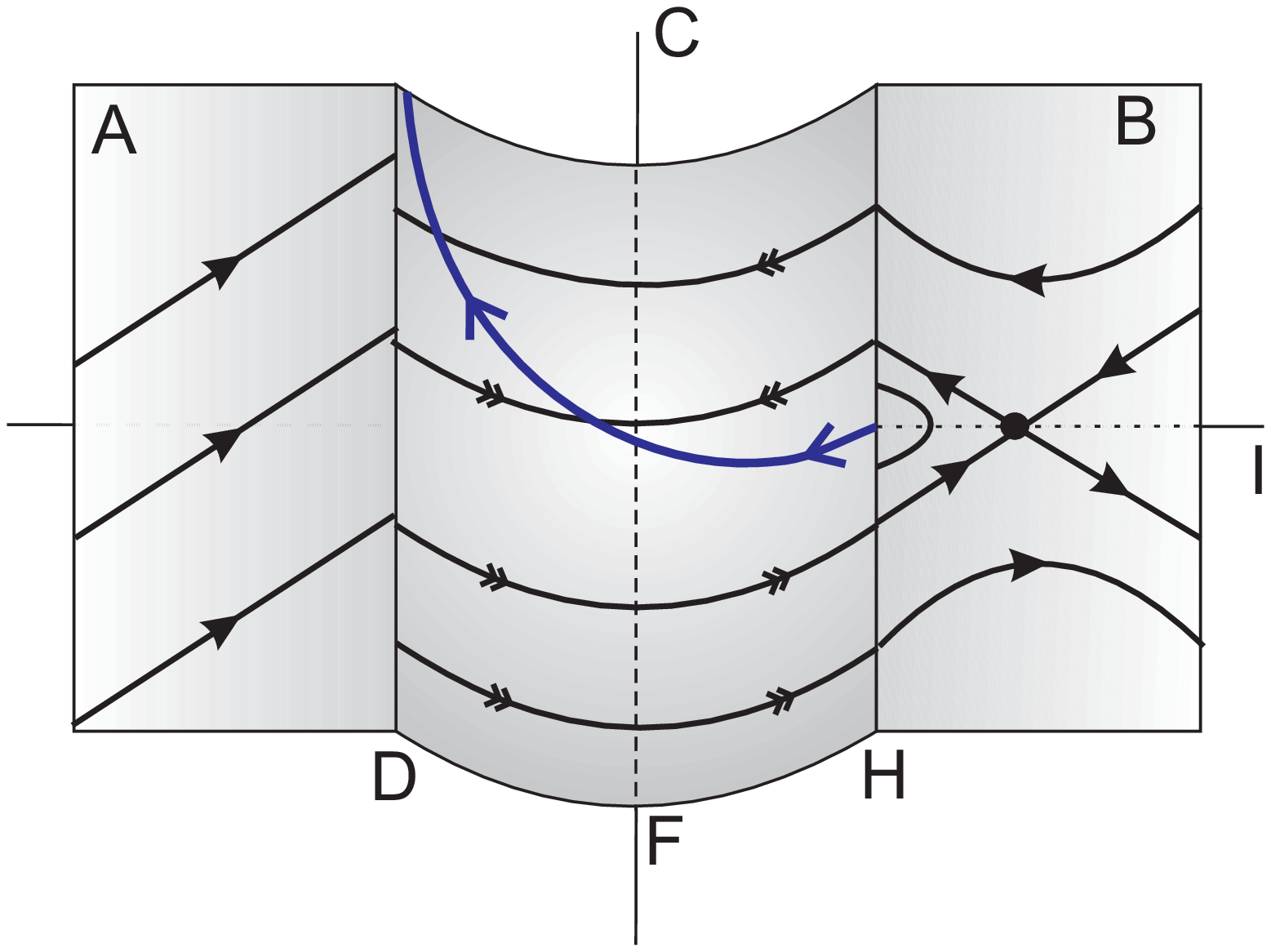}
\end{minipage}\caption{\small{In this figure is pictured the case $\lambda=0$ to the
left and the case $\lambda>0$ at the right. In both we consider
$\rho \in (0,\pi)$.}} \label{fig pert singular sela sobre a
descontinuidade e visivel}
\end{figure}

In Figure \ref{fig pert singular sem nada e invisivel} and in the
next ones,  double arrow over a curve 
means that it is a trajectory of the fast dynamical system, and
simple arrow means that it is a trajectory  of the one
dimensional slow dynamical system. 
%

%

\subsection{Regular$-$focus Bifurcation}

Consider the regular$-$focus normal form given in Subsection \ref{subsecao formas normais}.  
The regularized vector field becomes

$$\begin{array}{lcl}
   \dot{x} & = & \dfrac{1+\lambda+x+y}{2} +
\varphi \left( \frac{x}{\epsilon} \right)\dfrac{1-\lambda-x-y}{2} \,
,
\\
   \dot{y}  & =  & \dfrac{1-\lambda-x+y}{2} + \varphi \left(
\frac{x}{\epsilon} \right)\dfrac{1+\lambda+x-y}{2}.\end{array}
$$

Similarly to the previous case, considering the polar blow$-$up
given in $(\ref{blow-up polar})$, where $\theta = (2 \rho + \pi)/4$,
we get

\[
\begin{array}{lcl}
      r \dot{\theta}  & = & - \sin \theta \left(   \dfrac{1+\lambda+y+ r \cos \theta}{2} +
\varphi ( \cot \theta )\dfrac{1-\lambda  - y - r \cos \theta}{2} \right) \, , \\
  \dot{y} &= & \dfrac{1 - \lambda+y-r \cos \theta}{2} + \varphi (
\cot \theta)\dfrac{1+\lambda -y+r \cos \theta}{2}.
\end{array}
\]

Putting $r=0$, the fast dynamics is determined by the system
 \[ \theta'   = \sin \theta \left(   \dfrac{-1-\lambda-y}{2} +
\varphi ( \cot \theta )\dfrac{-1+\lambda+y}{2} \right) \, , \quad
y'=0 \, ;
\] and the slow dynamics on the slow manifold is determined by the
reduced system
\[    \dfrac{-1-\lambda-y}{2} +
\varphi ( \cot \theta )\dfrac{-1+\lambda+y}{2}   = 0 \, , \quad
\dot{y}= \dfrac{1-\lambda+y}{2} + \varphi ( \cot
\theta)\dfrac{1+\lambda-y}{2}.
\]

\begin{figure}[!h]
\begin{minipage}{0.49\linewidth}
\psfrag{A}{$X$} \psfrag{B}{$Y$} \psfrag{C}{$y$}
\psfrag{L}{$-\lambda$} \psfrag{D}{$0$} \psfrag{E}{$\frac{\pi}{4}$}
\psfrag{F}{$\frac{\pi}{2}$} \psfrag{G}{$\frac{3
\pi}{4}$}\psfrag{H}{$\pi$}\psfrag{I}{$x$} \epsfxsize=6cm
\epsfbox{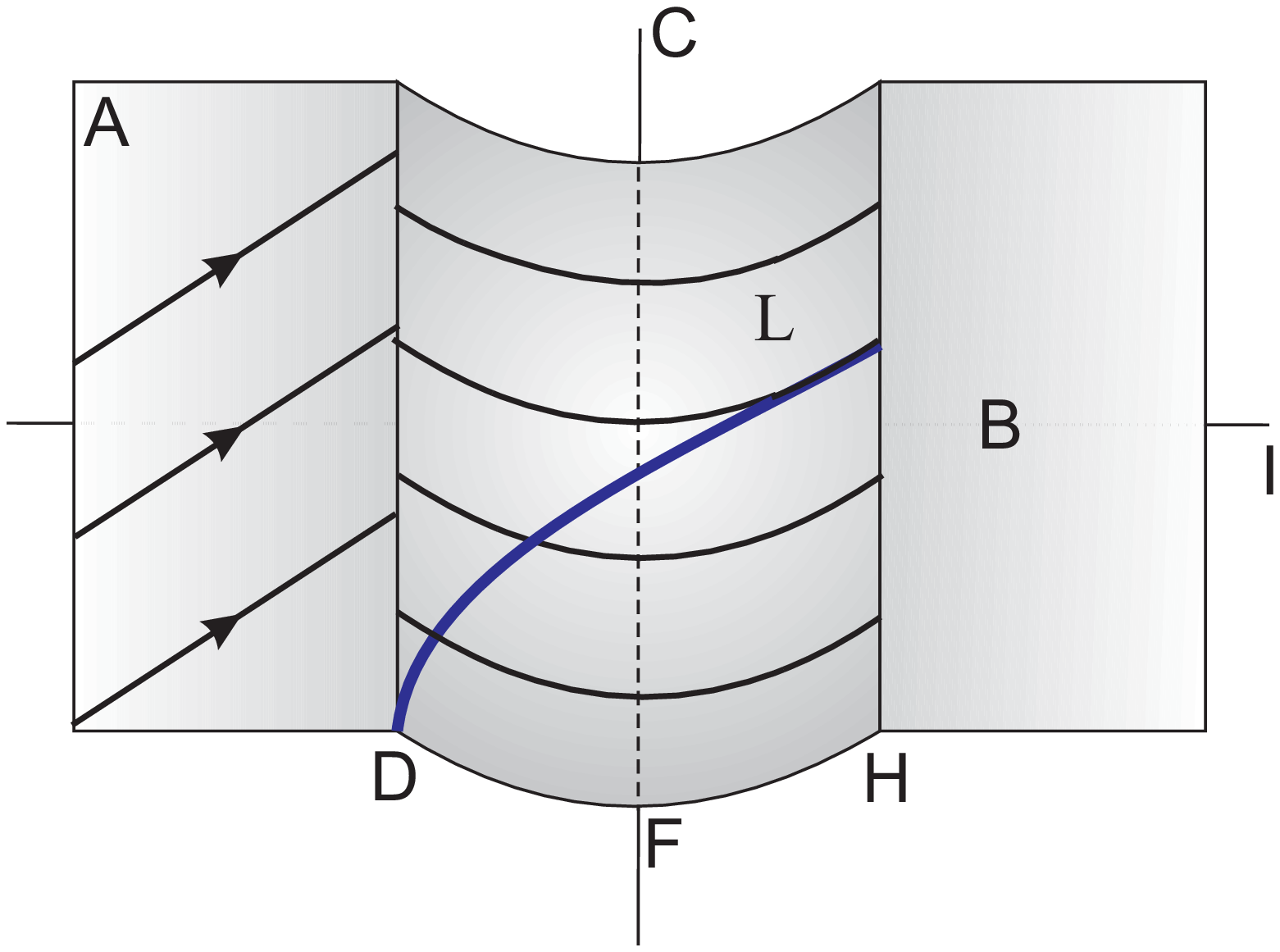}
\end{minipage} \hfill
\begin{minipage}{0.49\linewidth}
\psfrag{A}{$X$} \psfrag{B}{$Y$} \psfrag{C}{$y$}\psfrag{P}{$P$}
\psfrag{D}{$0$} \psfrag{E}{$\frac{\pi}{4}$}\psfrag{L}{$-\lambda$}
\psfrag{F}{$\frac{\pi}{2}$} \psfrag{G}{$\frac{3
\pi}{4}$}\psfrag{H}{$\pi$}\psfrag{I}{$x$} \epsfxsize=6.2cm
\epsfbox{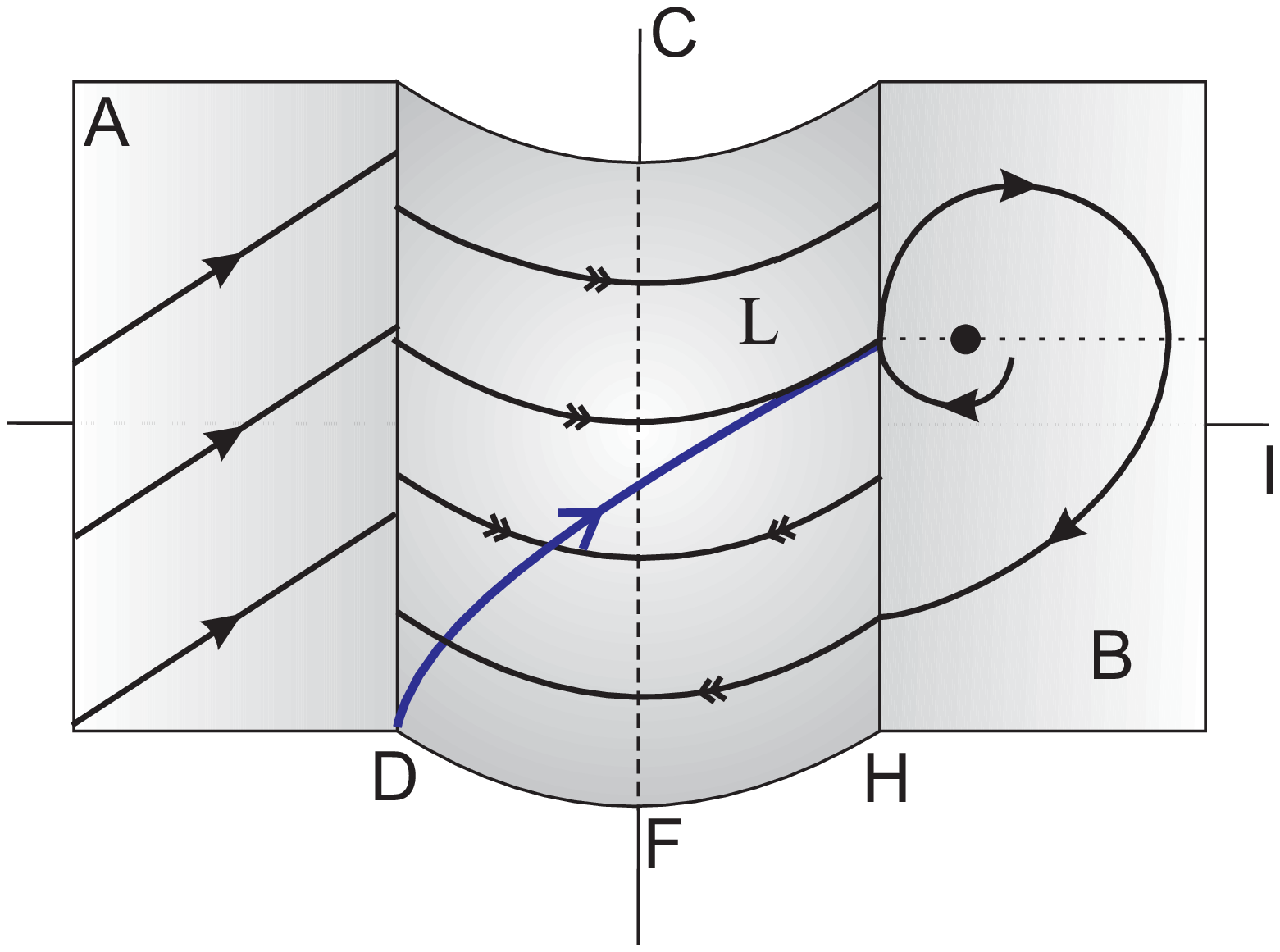}
\end{minipage}
\caption{\small{In this figure is pictured the slow manifold to the
left and the case $\lambda<0$ at the right. In both we consider
$\rho \in (0,\pi)$.}} \label{bif-cod1-regular-foco-vl-e-v3}
\end{figure}

In this case, the slow manifold depends of the parameter $\lambda$.
In fact, it is given implicitly by $(-1-\lambda-y)/2 + \varphi (\cot
\theta) (-1+\lambda+y)/2 = 0$ (see Figure
$\ref{bif-cod1-regular-foco-vl-e-v3}$). The slow manifold
$y(\theta)$, given in the previous equation, satisfies

\[
\displaystyle\lim_{\theta \longrightarrow \frac{\pi}{4}}y(\theta) =
- \infty \quad \mbox{and} \quad \displaystyle\lim_{\theta
\longrightarrow \frac{3 \pi}{4}}y(\theta) = - \lambda.
\]

\begin{figure}[!h]
\begin{minipage}{0.49\linewidth}
\psfrag{A}{$X$} \psfrag{B}{$Y$} \psfrag{C}{$y$}
\psfrag{L}{$-\lambda$} \psfrag{D}{$0$} \psfrag{E}{$\frac{\pi}{4}$}
\psfrag{F}{$\frac{\pi}{2}$} \psfrag{G}{$\frac{3
\pi}{4}$}\psfrag{H}{$\pi$}\psfrag{I}{$x$} \epsfxsize=6cm
\epsfbox{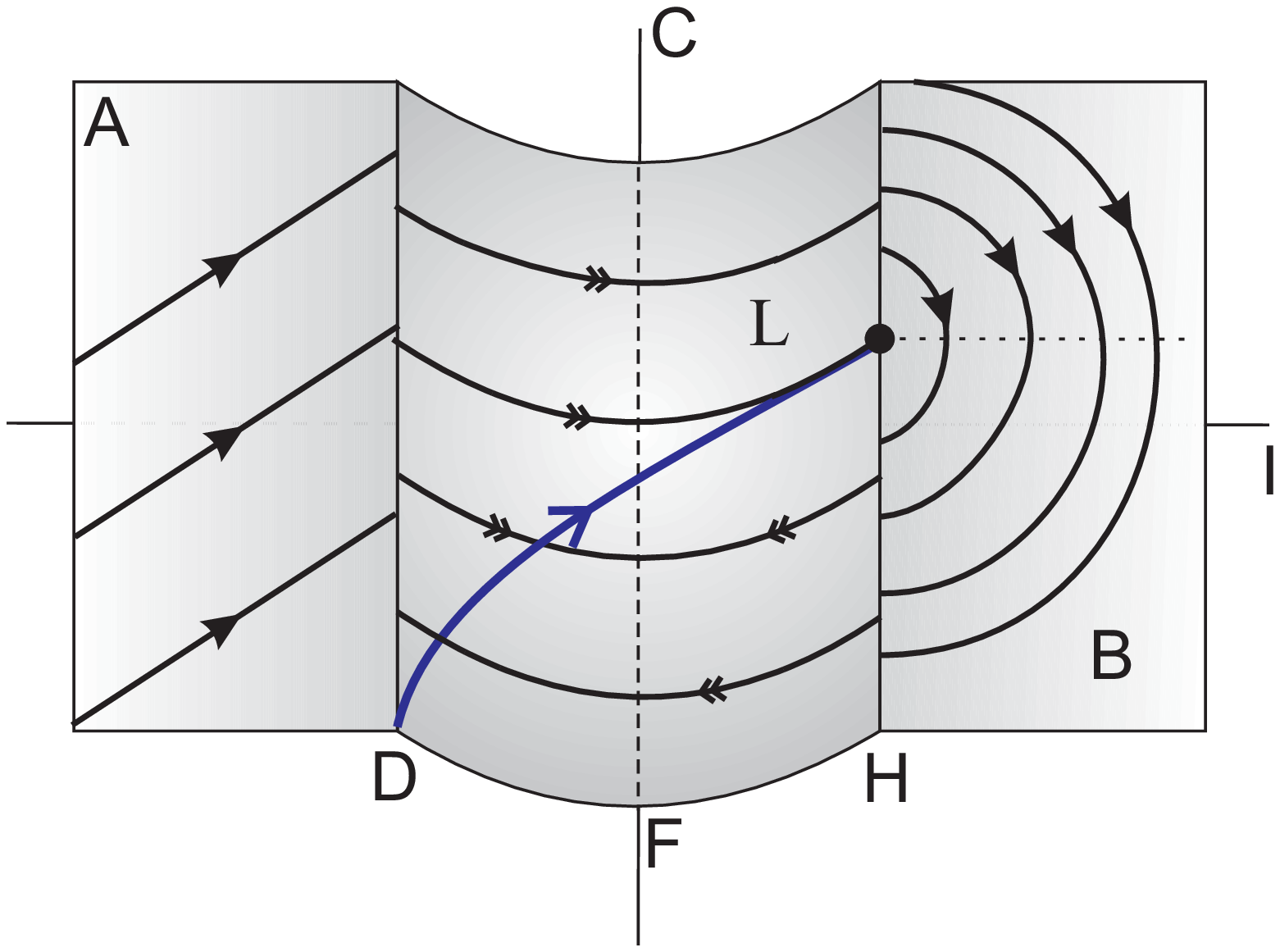}
\end{minipage} \hfill
\begin{minipage}{0.49\linewidth}
\psfrag{A}{$X$} \psfrag{B}{$Y$} \psfrag{C}{$y$}\psfrag{P}{$P$}
\psfrag{D}{$0$} \psfrag{E}{$\frac{\pi}{4}$}\psfrag{L}{$-\lambda$}
\psfrag{F}{$\frac{\pi}{2}$} \psfrag{G}{$\frac{3
\pi}{4}$}\psfrag{H}{$\pi$}\psfrag{I}{$x$} \epsfxsize=6.2cm
\epsfbox{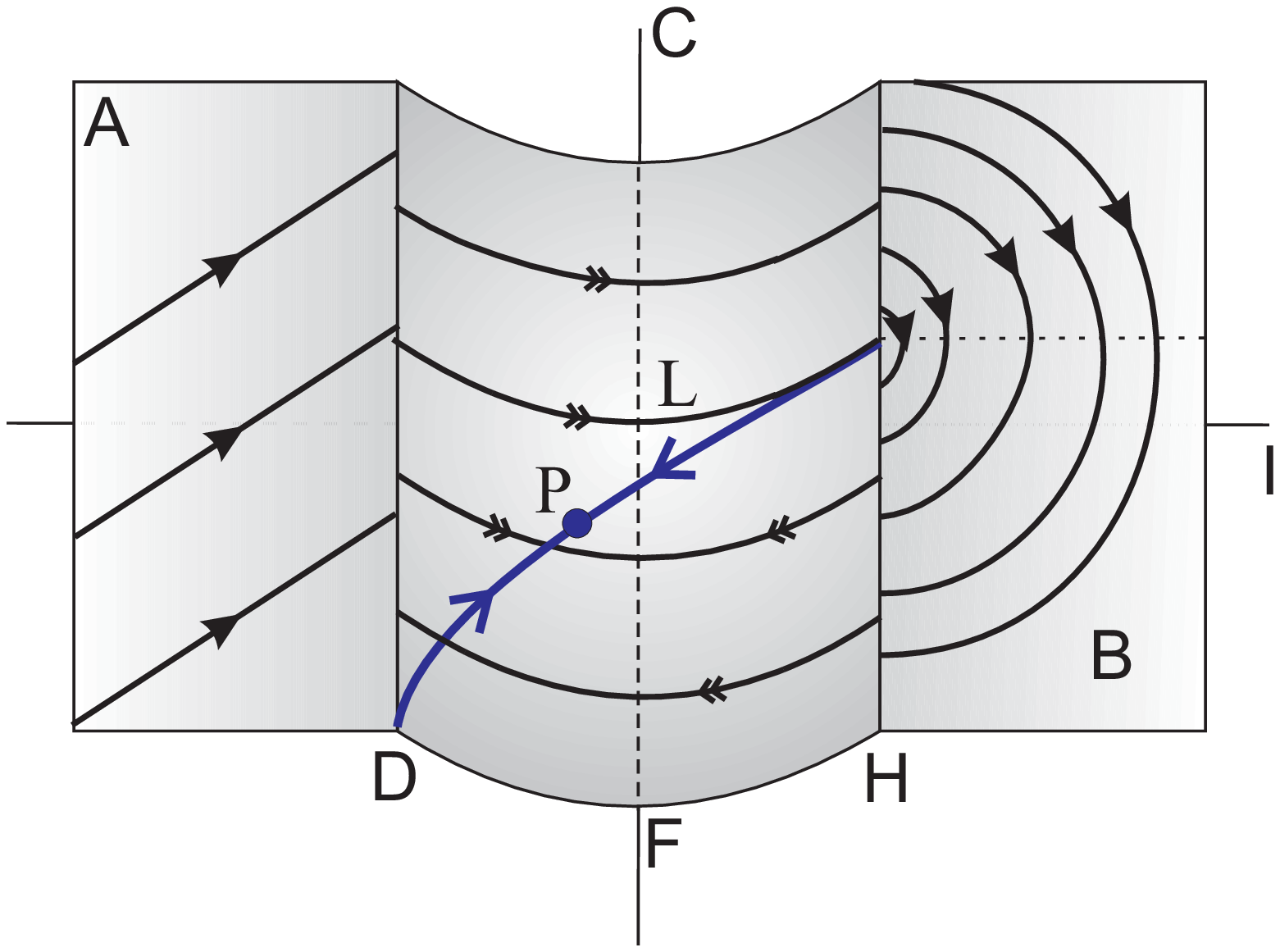}
\end{minipage}\caption{\small{In this figure is pictured the case $\lambda=0$ to the
left and the case $\lambda>0$ at the right. In both we consider
$\rho \in (0,\pi)$.}}\label{fig foco regular lambda zero e positivo}
\end{figure}

We give now the dynamics on the slow manifold. If $\lambda<0$ so
$\dot{y}\neq 0$ (see Figure \ref{bif-cod1-regular-foco-vl-e-v3}), if
$\lambda>0$ so $\dot{y}$ has an unique critical point $P$, given
implicitly as the solution of $\varphi ( \cot
\theta)=(-1+\lambda+y)/(1+\lambda-y)$,  which is an attractor (see
Figure \ref{fig foco regular lambda zero e positivo}) and if
$\lambda=0$ so $\dot{y}\neq 0$ (see Figure \ref{fig foco regular
lambda zero e positivo}).

\subsection{Regular$-$Cusp Bifurcation}

Consider the regular$-$cusp normal form given in Subsection \ref{subsecao formas normais}. 
The regularized vector field becomes

$$\begin{array}{lcl}
   \dot{x} & = & \dfrac{1+\lambda-y^2}{2}+\varphi \left( \dfrac{x}{\epsilon} \right)\dfrac{1-\lambda+y^2}{2}\,
,
\\\\
   \dot{y}  & =  & \varphi \left( \dfrac{x}{\epsilon} \right).\end{array}
$$

Making the polar blow$-$up given in $(\ref{blow-up polar})$, where
$\theta = (2 \rho + \pi)/4$, we get

\[
\begin{array}{lcl}
      r \dot{\theta}  & = & \dfrac{\sin\theta}{2} \, \left(\varphi(\cot \theta)(-1+\lambda-y^2)-1-\lambda+y^2\right) \, ,
      \\\\
  \dot{y} &= & \varphi(\cot \theta).
\end{array}
\]

Putting $r=0$ the fast dynamics is determined by the system

\[
\theta'   =
\dfrac{\sin\theta}{2}(\varphi(\cot\theta)(-1+\lambda-y^2)-1-\lambda+y^2)\,
, \quad y'=0 \, ;
\] and the slow dynamics on the slow manifold is determined by the
reduced system

\[
\varphi(\cot \theta)(-1+\lambda-y^2)-1-\lambda+y^2= 0 \, , \quad
\dot{y}= \varphi(\cot\theta).
\]

\begin{figure}[ht]
\epsfysize=10cm
\psfrag{A}{$X$}\psfrag{T}{$\rho^*$}\psfrag{P}{$p$}\psfrag{1}{$\lambda<0$}\psfrag{2}{$\lambda=0$}\psfrag{3}{$\lambda>0$}\psfrag{B}{$Y$}
\psfrag{Lp}{$\sqrt{\lambda}$}\psfrag{Ln}{$-\sqrt{\lambda}$}\psfrag{C}{$y$}
\psfrag{L}{$-\lambda$} \psfrag{D}{$0$} \psfrag{E}{$\frac{\pi}{4}$}
\psfrag{F}{$\frac{\pi}{2}$} \psfrag{G}{$\frac{3
\pi}{4}$}\psfrag{H}{$\pi$}\psfrag{I}{$x$}
\epsfbox{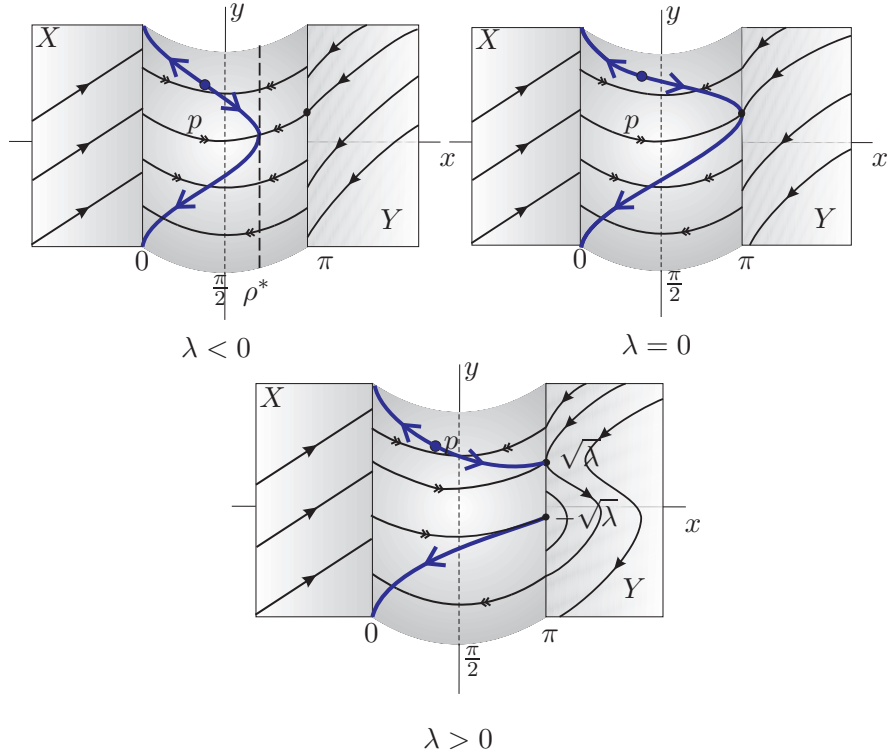}
\caption{\small{Regular$-$cusp bifurcation Diagram.}}
\label{bif-cod1-regular-cupide-vl-neg}
\end{figure}

Observe that the slow manifold depends of the parameter $\lambda$.
We can obtain the explicit form. In fact, the slow manifold is
composed by two branches (see Figure
$\ref{bif-cod1-regular-cupide-vl-neg}$):

\begin{equation}
y_{\pm}^{\lambda}(\theta)=\pm
\sqrt{\dfrac{\lambda(1-\varphi(\cot\theta))+1+\varphi(\cot\theta)}{1-\varphi(\cot\theta)}}.
\label{equacao-cuspide-regular}\end{equation}

The slow manifold satisfies the properties:

\begin{itemize}
\item [$(i)$] $\displaystyle\lim_{\theta \rightarrow \frac{\pi}{4}}
y_{\pm}^{\lambda}(\theta) = \pm \infty$;

\item [$(ii)$] If $\lambda<0$ there exists $\theta^*\in (\pi/4, 3\pi/4)$
(respectively, $\rho^* \in (0,\pi)$) such that
$y_{\pm}^{\lambda}(\theta^*)=0$ and the slow manifold is not defined
for $\rho \in (\rho^*,\pi)$. For $\theta\in (\pi/4, \theta^*)$ there
exist homeomorphisms $\xi_{\pm}$ between each branch of the slow
manifold and $\R^*$. That is, for each $z\in \R^*$ there exists
$\theta(z)\in (\pi/4,\theta^*)$ such that
$y_{\pm}^{\lambda}(\theta(z))=z$;

\item [$(iii)$] If $\lambda\geq 0$ there exist  homeomorphisms $\xi_{\pm}$
between each branch of the slow manifold and $\R^*$. That is, for
each $z\in \R^*$ there exists $\theta(z)\in (\pi/4, 3\pi/4)$ such
that $y_{\pm}^{\lambda}(\theta(z))=z$.

\end{itemize}

In fact, the item $(i)$ is a straightforward calculus. In order to
prove the item $(ii)$ observe Expression
$(\ref{equacao-cuspide-regular})$. Let  $\theta^*$ be such that
$\varphi(\cot \theta^*)=(1+\lambda)(\lambda-1)$. We have,

\[
y_{\pm}^{\lambda}(\theta^*)=0
\]and the radical in $(\ref{equacao-cuspide-regular})$ is negative
for $\theta\in (\theta^*,\pi)$.

We define the maps:

\begin{equation}
\begin{array}{lll}     \xi_{\pm}:       & \R^*      &\rightarrow (\pi/4,\theta^*) \\
                                  &  z       &\mapsto
                                  \theta(z)=\cot^{-1}\left(\varphi^{-1}\left(\dfrac{z^2-\lambda-1}{1-\lambda+_z^2}\right)\right).
\end{array}
\label{equacao-cuspide-regular-2}\end{equation}

Given $z\in \R^*$ if we put $\xi(z)=\theta(z)$ in
$(\ref{equacao-cuspide-regular-2})$ we get
$y_{\pm}^{\lambda}(\theta(z))=z$. Note that $\xi_{\pm}$ are
homeomorphisms.

The proof of the item $(iii)$ is analogous.

The dynamics in the slow manifold is given by
$\dot{y}=\varphi(\cot\theta)$. So there exists a unique critical
point $p$ given implicitly as the solution of $\varphi(\cot
\theta_p) = 0$. Note that this critical point is a repeller because
$\varphi(\cot \theta) < 0$ for $\theta < \theta_p$ and $\varphi(\cot
\theta) > 0$ for $\theta>\theta_p$. See Figure
\ref{bif-cod1-regular-cupide-vl-neg}.
%

\section{Fold$-$fold Bifurcations}\label{secao fold fold}

\begin{figure}[ht]
\epsfysize=3.5cm
\psfrag{A}{$X$}\psfrag{B}{$Y$}\psfrag{X}{$x$}\psfrag{Y}{$y$}\psfrag{1}{$(a)$}\psfrag{2}{$(b)$}\psfrag{3}{$(c.1)$}\psfrag{4}{$(c.2)$}
\epsfbox{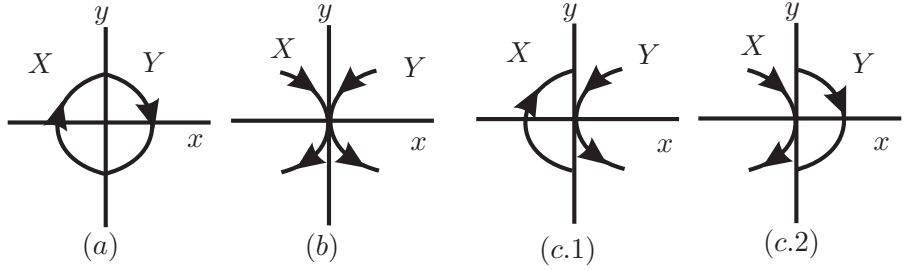} \caption{\small{Fold$-$fold
singularities.}}\label{dobras}
\end{figure}

In this section we analyze the dynamics of a NSDS around a point $q$
which is $\Sigma-$fold of both $X$ and $Y$. We say that $q$ is a
\textit{Fold$-$Fold singularity} of $Z \in \Omega^r$. We divide the
fold$-$fold singularities in types according with the sign of
$X^2.f(q)$ and $Y^2.f(q)$:

\begin{itemize}
\item [$(a)$] \textbf{Elliptic case}: $X^2.f(q)>0$ and $Y^2.f(q)<0$. 
See Figure $\ref{dobras}$ (a).

\item [$(b)$] \textbf{Hyperbolic case}: $X^2.f(q)<0, Y^2.f(q)>0$. See Figure $\ref{dobras}$ (b).

\item [$(c.1)$] \textbf{Parabolic visible case}: $X^2.f(q)>0,
Y^2.f(q)>0$. See Figure $\ref{dobras}$ (c.1).

\item [$(c.2)$] \textbf{Parabolic invisible case}:
$X^2.f(q)<0,Y^2.f(q)<0$. See Figure $\ref{dobras}$ (c.2).
\end{itemize}

Note that, we can define a first return map $\psi_Z$ only in
the elliptic case. 
Take $\Sigma$ as the $y-$axis, i.e., $f(x,y)=x$ and consider the
parameter $\lambda \in (-1,1)$. The generic normal forms of the
hyperbolic and parabolic fold$-$fold singularities are given in
Table \ref{tabela}. The normal form of the elliptic fold$-$fold
singularity is given in Subsection \ref{fold-fold elliptic}.

\begin{table}[ht!]
\begin{tabular}{|c|c|c|}              \hline
                         \textbf{Hyperbolic}                    &\textbf{Parabolic visible}    &\textbf{Parabolic invisible}
                           \\  \hline \hline
                   $X_{\lambda}(x,y)=(y-\lambda,-1)$     & $X_{\lambda}(x,y)=(y-\lambda,1)$    &$X_{\lambda}(x,y)=(y-\lambda,-1)$\\
                   $Y_{}(x,y)=(-y,-1)$            & $Y_{}(x,y)=(-y,-1)$          &$Y_{}(x,y)=(-y,1)$  \\
                                     \hline
\end{tabular}\caption{}\label{tabela}
\end{table}

%
In the next three subsections we study the dynamics of the
hyperbolic and parabolic fold$-$fold singularities via geometric
singular perturbations.

\subsection{Hyperbolic Case}

Consider the normal form of the hyperbolic fold$-$fold singularity
given in Table \ref{tabela}. The regularized vector field is given
by

$$\begin{array}{lcl}
   \dot{x}            &= & -\dfrac{\lambda}{2}+ \varphi\left(\dfrac{x}{\epsilon} \right)\left(\dfrac{-\lambda+2y}{2}\right) \, ,
\\\\
   \dot{y}  & =  & -1.\end{array}
$$ By the polar blow up we get

\[
\begin{array}{lcl}  r\dot{\theta}  & = & \sin \theta \left( \dfrac{\lambda}{2} +
\varphi(\cot \theta) \dfrac{\lambda -2y}{2} \right) \, , \\
  \dot{y} &= & -1
\end{array}
\] where $\theta = (2 \rho + \pi)/4$.

Putting $r=0$ the fast dynamics is determined by the system

 \[
 \theta'   = \sin \theta \left( \dfrac{\lambda}{2} +
\varphi(\cot \theta) \left(\dfrac{\lambda -2y}{2}\right) \right) \,
, \quad y'=0 \, ;
\] and the slow dynamics on the slow manifold is determined by the
reduced system

\[
\dfrac{\lambda}{2} + \varphi(\cot \theta) \left(\dfrac{\lambda
-2y}{2}\right)=0 \, , \quad \dot{y}= -1.
\]

In this case we obtain the explicit expression for the slow
manifold:

\[
y(\theta)=\frac{\lambda(1+\varphi(\cot
\theta))}{2\varphi(\cot\theta)}.
\]Observe that, the slow manifold $y(\theta)$ is not
defined for $\theta_0$ such that $\varphi(\cot \theta_0)=0$. So, for
$\lambda\neq 0$, $y(\theta)$ have two branches and satisfies:

\begin{itemize}


    \item [$(a)$] $\displaystyle\lim_{\theta \longrightarrow \theta_{0}^-}
    y(\theta)=-\infty$ for $\lambda<0$ and $\displaystyle\lim_{\theta \longrightarrow
    \theta_{0}^-} y(\theta)=+\infty$ for $\lambda>0$;

    \item [$(b)$] $\displaystyle\lim_{\theta \longrightarrow \theta_{0}^+}
    y(\theta)=+\infty$ for $\lambda<0$ and $\displaystyle\lim_{\theta \longrightarrow
    \theta_{0}^+} y(\theta)=-\infty$ for $\lambda>0$.

    \item[$(c)$] For $\lambda=0$ the slow manifold is given implicitly by $y \varphi(\cot
    \theta)=0$, that is, $\{ (\theta, y) \, | \, \theta=\theta_0 \} \cup \{ (\theta, y) \, | \, y=0 \}$ is the slow
    manifold.
\end{itemize}

\begin{figure}[ht]
\epsfysize=10cm
\psfrag{0}{$0$}\psfrag{A}{$X_{\lambda}$}\psfrag{T}{$\theta_0$}\psfrag{P}{$p$}\psfrag{1}{$\lambda<0$}\psfrag{2}{$\lambda=0$}\psfrag{3}{$\lambda>0$}\psfrag{B}{$Y$}
\psfrag{Lp}{$\sqrt{\lambda}$}\psfrag{Ln}{$-\sqrt{\lambda}$}\psfrag{C}{$y$}
\psfrag{L}{$\lambda$} \psfrag{D}{$0$} \psfrag{E}{$\frac{\pi}{4}$}
\psfrag{F}{$\frac{\pi}{2}$} \psfrag{G}{$\frac{3
\pi}{4}$}\psfrag{H}{$\pi$}\psfrag{I}{$x$}
\epsfbox{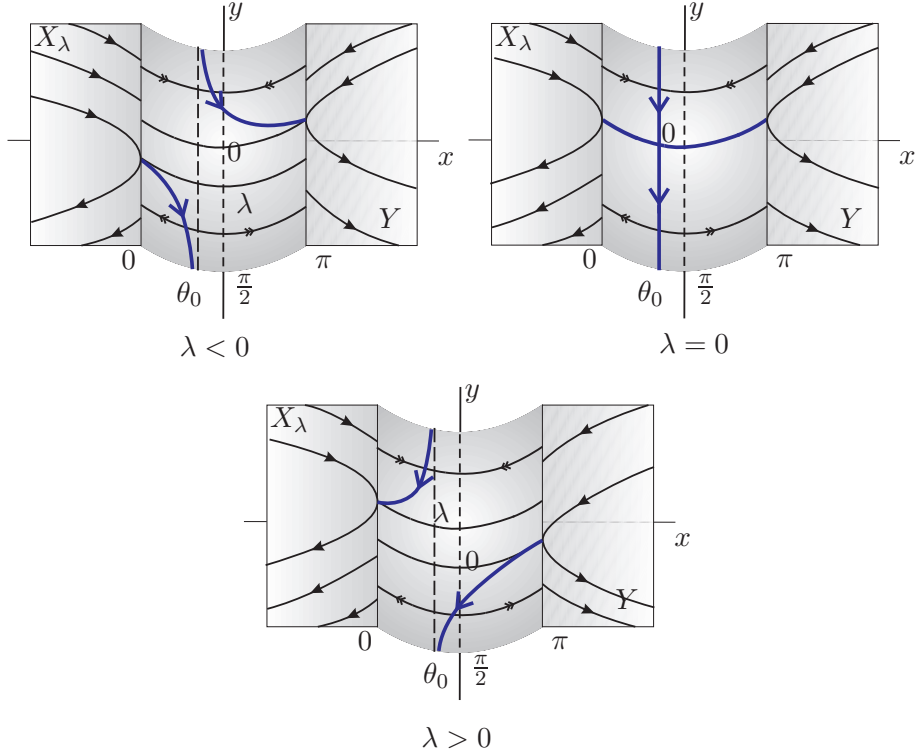} \caption{\small{Slow
manifold depending of the parameter $\lambda$.}}
\label{bif-cod1-fold-hiperbolico-vl}
\end{figure}

The dynamics on the slow manifold is given by $\dot{y}=-1$.
Therefore, do not exist critical points. See Figure
$\ref{bif-cod1-fold-hiperbolico-vl}$.

\subsection{Parabolic visible case}

Consider the normal form of the parabolic visible fold$-$fold
singularity given in Table \ref{tabela}.
The regularized vector field is

\[
\begin{array}{lcl}
   \dot{x}            &= & -\dfrac{\lambda}{2}+ \varphi\left(\dfrac{x}{\epsilon} \right)
   \left(\dfrac{-\lambda+2y}{2}\right) \, ,
\\\\
   \dot{y}  & =  & \varphi\left(\dfrac{x}{\epsilon} \right).\end{array}
\] By the polar blow up we get

\[
\begin{array}{lcl}  r\dot{\theta}  & = & \sin \theta \left( \dfrac{\lambda}{2} +
\varphi(\cot \theta)\left( \dfrac{\lambda -2y}{2}\right) \right) \,
,
\\\\
  \dot{y} &= & \varphi(\cot \theta).
\end{array}
\]where $\theta = (2 \rho + \pi)/4$.

Putting $r=0$ the fast dynamics is determined by the system

 \[
 \theta'   = \sin \theta \left( \dfrac{\lambda}{2} +
\varphi(\cot \theta) \left(\dfrac{\lambda -2y}{2}\right) \right) \,
, \quad y'=0 \, ;
\] and the slow dynamics on the slow manifold is determined by the
reduced system

\[
\dfrac{\lambda}{2} + \varphi(\cot \theta) \left(\dfrac{\lambda
-2y}{2}\right)=0 \, , \quad \dot{y}= \varphi(\cot \theta).
\]

\begin{figure}[ht]
\epsfysize=10cm
\psfrag{0}{$0$}\psfrag{A}{$X$}\psfrag{T}{$\theta_{0}$}
\psfrag{P}{$p$}\psfrag{1}{$\lambda<0$}\psfrag{2}{$\lambda=0$}\psfrag{3}{$\lambda>0$}\psfrag{B}{$Y$}
\psfrag{Lp}{$\sqrt{\lambda}$}\psfrag{Ln}{$-\sqrt{\lambda}$}\psfrag{C}{$y$}
\psfrag{L}{$\lambda$} \psfrag{D}{$0$} \psfrag{E}{$\frac{\pi}{4}$}
\psfrag{F}{$\frac{\pi}{2}$} \psfrag{G}{$\frac{3
\pi}{4}$}\psfrag{H}{$\pi$}\psfrag{I}{$x$}
\epsfbox{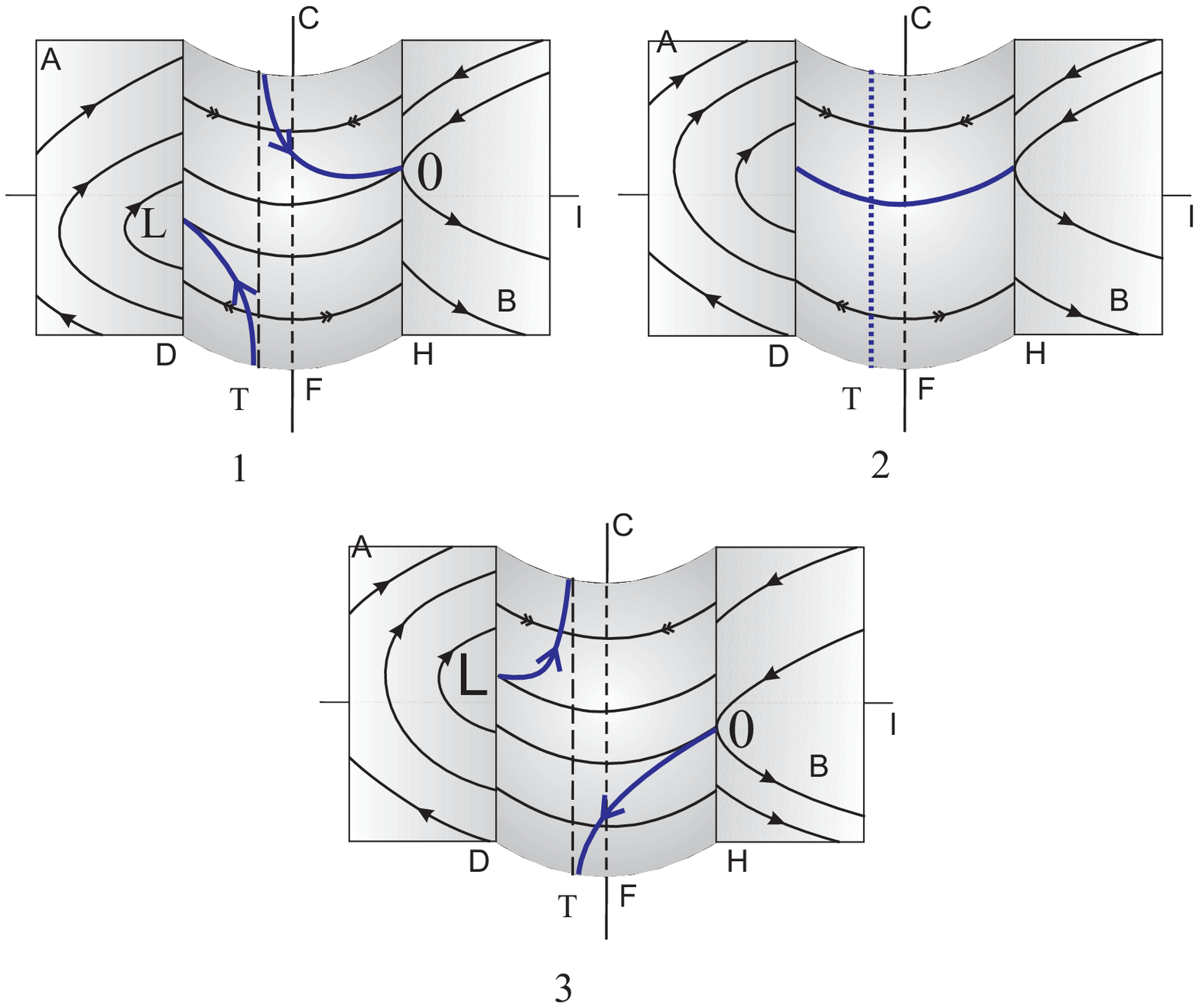}
\caption{\small{Bifurcation Diagram of the Parabolic Visible
Fold$-$Fold Singularity.}}
\label{bif-cod1-fold-parabolico-visivel-vl}
\end{figure}

The analysis is similar to the hyperbolic case. In the present case
the dynamics on the slow manifold is given by $\dot{y}=\varphi(\cot
\theta )$. So, for $\lambda=0$, the straight line
$\theta=\theta_{0}$ is composed by critical points. See Figure
$\ref{bif-cod1-fold-parabolico-visivel-vl}$.

\subsection{Parabolic invisible case}

For this case, we get one different topo\-lo\-gi\-cal type of
bifurcation. Consider the normal form of the parabolic invisible fold$-$fold singularity given in Table \ref{tabela}. 
The regularized vector field is

\[
\begin{array}{lcl}
   \dot{x}            &= & -\dfrac{\lambda}{2} \varphi\left(\dfrac{x}{\epsilon} \right)+\dfrac{-\lambda+2y}{2} \, ,
\\\\
   \dot{y}  & =  & -1 .\end{array}
\] By the polar blow up we get

\[
\begin{array}{lcl}  r\dot{\theta}  & = & \dfrac{\sin \theta}{2} \left(\lambda\varphi(\cot \theta)+\lambda -2y\right) \, ,
\\\\
  \dot{y} &= & -1.
\end{array}
\]where $\theta = (2 \rho + \pi)/4$.

Putting $r=0$ the fast dynamics is determined by the system

 \[
 \theta'   = \dfrac{\sin \theta}{2} \left(\lambda\varphi(\cot \theta)+\lambda -2y\right) \,
, \quad y'=0 \, ;
\] and the slow dynamics on the slow manifold is determined by the
reduced system

\[
\sin \theta \left( \dfrac{\lambda}{2}\varphi(\cot \theta)+\left(
\dfrac{\lambda -2y}{2}\right) \right)=0 \, , \quad \dot{y}=-1.
\]

We have the explicit expression for the slow manifold in this case:

\[
y(\theta)=\dfrac{\lambda}{2}(1+\varphi(\cot \theta)).
\]

\begin{figure}[ht]
\epsfysize=10cm
\psfrag{0}{$0$}\psfrag{A}{$X_{\lambda}$}\psfrag{T}{$\widetilde{\theta}$}
\psfrag{P}{$p$}\psfrag{1}{$\lambda<0$}\psfrag{2}{$\lambda=0$}\psfrag{3}{$\lambda>0$}
\psfrag{B}{$Y_{\lambda}$}\psfrag{Lp}{$\sqrt{\lambda}$}\psfrag{Ln}{$-\sqrt{\lambda}$}\psfrag{C}{$y$}
\psfrag{L}{$\lambda$} \psfrag{D}{$0$} \psfrag{E}{$\frac{\pi}{4}$}
\psfrag{F}{$\frac{\pi}{2}$} \psfrag{G}{$\frac{3
\pi}{4}$}\psfrag{H}{$\pi$}\psfrag{I}{$x$}
\epsfbox{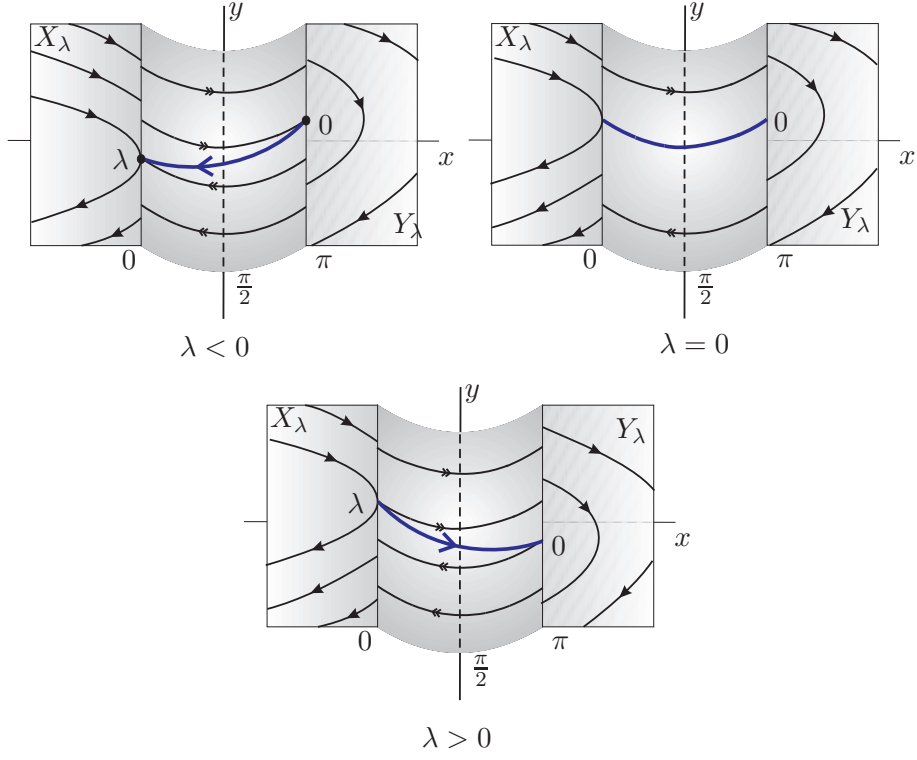}
\caption{\small{Bifurcation Diagram of the Parabolic Invisible
Fold$-$Fold Singularity.}}
\label{bif-cod1-fold-parabolico-invisivel-vl}
\end{figure}

The analysis is similar to the previous cases and the bifurcation
diagram is expressed in Figure
$\ref{bif-cod1-fold-parabolico-invisivel-vl}$.

\subsection{Elliptic case}\label{fold-fold elliptic}

In this case, associated with the non$-$smooth vector fields, there
exist the first return map $\psi_Z(p)$. Therefore, we need to
analyze the structural stability of this one dimensional
diffeomorphism.

Consider $Z$ presenting an elliptic fold$-$fold
singularity, $f(x,y)=x$ and 

\begin{equation}\label{desdobramento-ff-eliptico}
Z_{\lambda}(x,y)=\left\{
                   \begin{array}{ll}
                     X_{\lambda}(x,y)=(y-\lambda,1), & \hbox{ for $(x,y) \in \Sigma_+$,} \\
                     Y(x,y)=(y,-1), & \hbox{ for $(x,y) \in \Sigma_-$.}
                   \end{array}
                 \right.
\end{equation} The regularized
vector field is

\[
\begin{array}{lcl}
   \dot{x}            &= & -\varphi\left(\dfrac{x}{\epsilon} \right)\dfrac{\lambda }{2}+
   \dfrac{2y-\lambda}{2} \, ,
\\\\
   \dot{y}  & =  & \varphi\left(\dfrac{x}{\epsilon} \right).\end{array}
\] By the polar blow up we get

\[
\begin{array}{lcl}  r\dot{\theta}  & = & \sin \theta \left( \dfrac{\lambda \varphi(\cot \theta)}{2}+
\left( \dfrac{\lambda -2y}{2}\right) \right) \, , \\\\
  \dot{y} &= & \varphi(\cot \theta)
\end{array}
\]where $\theta = (2 \rho + \pi)/4$.

Putting $r=0$ the fast dynamics is determined by the system

\[
 \theta' = \sin \theta \left( \dfrac{\lambda \varphi(\cot \theta)}{2}+
 \left(\dfrac{\lambda -2y}{2}\right) \right) \, , \quad y'=0 \, ;
\] and the slow dynamics on the slow manifold is determined by the
reduced system

\[
\dfrac{\lambda\varphi(\cot \theta)}{2} + \left(\dfrac{\lambda
-2y}{2}\right)=0 \, , \quad \dot{y}= \varphi(\cot \theta).
\]

\begin{figure}[ht]
\epsfysize=10cm
\psfrag{0}{$0$}\psfrag{A}{$X_{\lambda}$}\psfrag{T}{$\widetilde{\theta}$}
\psfrag{P}{$p$}\psfrag{1}{$\lambda<0$}\psfrag{2}{$\lambda=0$}\psfrag{3}{$\lambda>0$}
\psfrag{B}{$Y_{\lambda}$}\psfrag{Lp}{$\sqrt{\lambda}$}\psfrag{Ln}{$-\sqrt{\lambda}$}\psfrag{C}{$y$}
\psfrag{L}{$\lambda$} \psfrag{D}{$0$} \psfrag{E}{$\frac{\pi}{4}$}
\psfrag{F}{$\frac{\pi}{2}$} \psfrag{G}{$\frac{3
\pi}{4}$}\psfrag{H}{$\pi$}\psfrag{I}{$x$}
\epsfbox{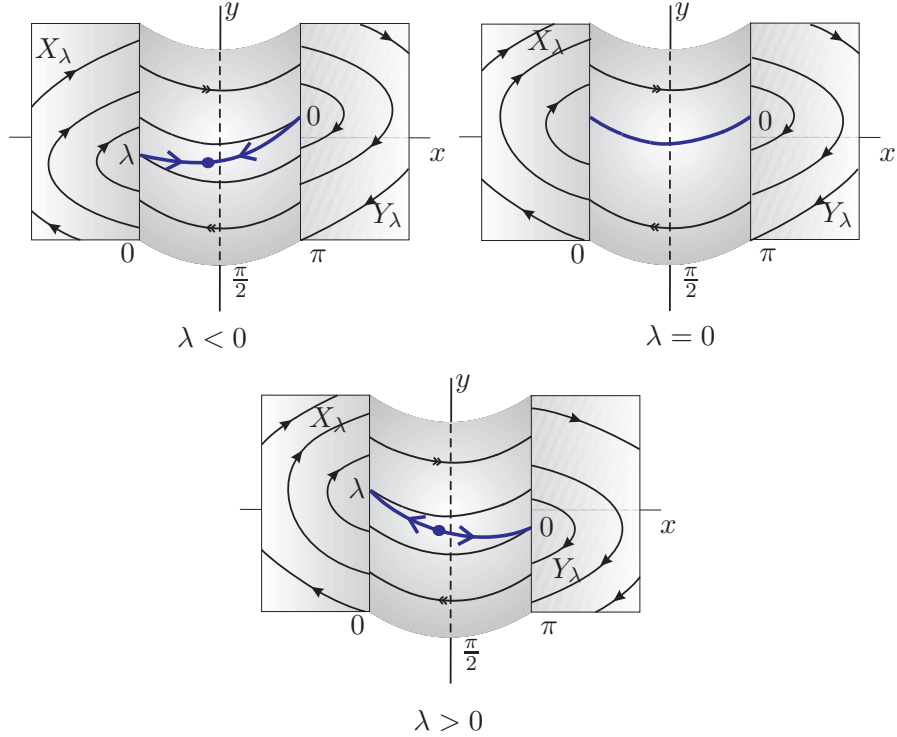}
\caption{\small{Bifurcation Diagram of the $1-$parameter Elliptic
Fold$-$Fold Singularity.}} \label{bif-cod1-fold-eliptico-vl}
\end{figure}

In this case, for $\lambda=0$, we only have sewing region on the
non$-$smooth manifold $\Sigma$. The explicit expression for the slow
manifold is
\[
y(\theta)=\dfrac{\lambda(1+\varphi(\cot \theta))}{2}
\]and there exist one critical point which is attractor if
$\lambda<0$ and repeller if $\lambda>0$. See Figure
$\ref{bif-cod1-fold-eliptico-vl}$.

Differently of the hyperbolic and parabolic cases, the unfolding
$(\ref{desdobramento-ff-eliptico})$ does not give the generic
unfolding of a non$-$smooth vector field presenting a elliptic
fold$-$fold singularity. 

%
%
%
%
%
%
%
%
%

Let $Z_{\lambda}=(X_{\lambda},Y)$ be as in
$(\ref{desdobramento-ff-eliptico})$. The expression of its first
return map is

\begin{equation}
\psi_{Z_{\lambda}}(y)=\gamma_Y\circ
\gamma_{X_{\lambda}}(y)=y-2\lambda,
\label{expressao-pr}\end{equation}where $\gamma_{X_{\lambda}}(p)$
(respectively $\gamma_Y(p)$) is the first return to $\Sigma$ of the
trajectory of $X_{\lambda}$ (respectively $Y$) that passes through
$p$.

Therefore, 
we conclude that the critical point of $\psi_{Z_{\lambda}}$ is not
hyperbolic. In order to obtain the generic unfolding of this case we
need to  unfold  the first return map. So, the unfolding of the
elliptic fold$-$fold singularity depends on two parameters. The
first, $\lambda$, is responsible by the displacement of one fold
along the $y-$axis and another one, $\varepsilon$, for the unfolding
of $\psi_{Z_{\lambda}}$.

Consider the $\varepsilon-$perturbation of the smooth vector field
$Y$, given by:

\begin{equation}
\begin{array}{ll}   Y_{\varepsilon}(x,y)            &=g_1^{\varepsilon}(x(t),
y(t))\dfrac{\partial}{\partial x}+g_2^{\varepsilon}(x(t),
y(t))\dfrac{\partial}{\partial y}\\\\
                                                   & =g_1^{\varepsilon}(t)\dfrac{\partial}{\partial
x}+g_2^{\varepsilon}(t)\dfrac{\partial}{\partial y}.
\end{array}
\label{desdobramento-campoY}\end{equation}

So, the flows of $X_{\lambda}$ and $Y_{\varepsilon}$ are:

\[
\begin{array}{ll} \phi_{X_{\lambda}}^t(x_0,y_0)             &=(x_0+(y_0-\lambda)t+t^2/2,
y_0+t),\\\\
               \phi_{Y_{\varepsilon}}^t(x_0,y_0)            &=\left(x_0+\int_0^tg_1^{\varepsilon}(s)\,ds,
y_0+\int_0^tg_2^{\varepsilon}(s)\,ds \right).
\end{array}
\]

Let $t^*\in \R^*$ and $t_1=2(\lambda-y_0)$ such that

\begin{equation}
\int_0^{t^*}g_1^{\varepsilon}(s)\,ds=0
\label{equacao-fluxo-Y}\end{equation}and
$\phi_{X_{\lambda}}^{t_1}(0,y_0)=(0, -y_0+2\lambda)\in \Sigma$.

Observe that there exist $t^*$ as in Equation
$(\ref{equacao-fluxo-Y})$, because $0$ is a elliptic fold
singularity for $Y$.

We suppose that $g_i^{\varepsilon}(.), i=1,2$ satisfies:

\begin{itemize}
    \item [$(a)$] $g_i$ are $C^r$ functions for $i=1,2$;
    \item [$(b)$] $Y_{\varepsilon}.f(0,0)=g_1^{\varepsilon}(0)=0$;
    \item [$(c)$] $Y_{\varepsilon}^2.f(0,0)=g_1^{\varepsilon}(0)\dfrac{d}{dx} g_1^{\varepsilon}(0)
    +g_2^{\varepsilon}(0)\dfrac{d}{dy} g_1^{\varepsilon}(0)\neq 0$;
    \item [$(d)$] $\int_0^{t^*}g_2^{\varepsilon}(s)\,
    ds=(\varepsilon-2)y+O(y^2)$.
\end{itemize}

The smooth vector fields $X_{\lambda}, Y_{\varepsilon}$ exhibited in
$(\ref{desdobramento-ff-eliptico}), (\ref{desdobramento-campoY})$,
respectively, supply the unfolding of $\psi_Z$:

\[
\psi_{Z_{\lambda,\varepsilon}}(y)=\phi_{Y_{\varepsilon}}^{t^*}\circ
\phi_{X_{\lambda}}^{t_1}(0,y)=(1-\varepsilon)y-2\lambda+2\lambda
\varepsilon+O(y^2).
\]

Therefore, the generic unfolding for the non$-$smooth vector field
$Z$ with the origin is an elliptic fold$-$fold singularity is
$Z_{\lambda, \varepsilon}=(X_{\lambda}, Y_{\varepsilon})$ where

\begin{equation}\label{forma normal fold-fold eliptica}
Z_{\lambda, \varepsilon}(x,y)= \left\{
                                 \begin{array}{ll}
                                   X_{\lambda}(x,y)=(y-\lambda,1), & \hbox{ if $(x,y) \in \Sigma_{+}$,} \\
                                   Y_{\varepsilon}(x,y)=(g_1^{\varepsilon}(x,y),
g_2^{\varepsilon}(x,y)), & \hbox{ if $(x,y) \in \Sigma_{-}$}
                                 \end{array}
                               \right.
\end{equation}and the smooth function $g_i^{\varepsilon},i=1,2$ satisfies the
conditions $(a),(b),(c)$ and $(d)$ given previously.

\section{Conclusion}\label{secao conclusao}

We note that: if for any $q \in \Sigma$ we have that $Xf(q)\neq 0$
or $Yf(q) \neq 0$ then, by Theorem 1.1 of \cite{LST}, there exists a
singular perturbation problem such that the sliding region is
homeomorphic to the slow manifold and the sliding vector field is
topologically equivalent to the reduced problem. This fact is useful in the next two subsections.\\

\subsection{Proof of Theorem \ref{teorema boundary bifurcations}}
In face of Theorem 1.1 of \cite{LST}, this theorem is the subject of
section \ref{secao boundary bifurcations}. Moreover, as we give the
topological behavior of the cases $\lambda<0$, $\lambda=0$ and
$\lambda>0$ it is easy to construct the
bifurcation diagram of \ref{eq pert sing}.\\

\subsection{Proof of Theorem \ref{teorema bif fold-fold}}
In this theorem we extend Theorem 1.1 of \cite{LST} considering that
can exists a
point $q$ such that 
$X.f(q)=Y.f(q)=0$, $X^{2}.f(q)\neq 0$ and $Y^{2}.f(q)\neq0$. In this
way,  $q$ is a 
$\Sigma-$fold point of both $X$ and $Y$. 
%

Consider a NSDS $Z_\lambda=(X,Y)$ where $\lambda \in \R$ is a
parameter. If with the variation of $\lambda \in (- \varepsilon,
 \varepsilon)$, the following behaviors are observable
then we consider that $Z_{\lambda}$ presents a bifurcation, where
$\varepsilon>0$ and small. Consider $\lambda^{+} \in
(0,\varepsilon)$ and $\lambda^{-} \in (-\varepsilon,0)$. The
behaviors are:
\begin{description}
\item[(i)] A change of stability on $\Sigma$, i.e., where $Z_{\lambda^{+}}$ has
a sliding region $\Sigma_3$ the non$-$smooth vector field
$Z_{\lambda^{-}}$ has an escaping region $\Sigma_2$.

\item[(ii)] A change of stability on $\dot{y}_{\lambda}$, i.e., there are
components of $\Sigma$ such that the induced flow on the slow
manifold is such that $\dot{y}_{\lambda^{+}}>0$ and
$\dot{y}_{\lambda^{-}}<0$.

\item[(iii)] A change of stability of the $\Sigma-$singularity,
i.e.,  $Z_{\lambda^{+}}$ presents a $\Sigma-$attractor and
$Z_{\lambda^{-}}$ presents a $\Sigma-$repeller.

\item[(iv)] A change of orientation on $\Sigma_1$ (the sewer
region), i.e., $Z_{\lambda^{+}}$ and $Z_{\lambda^{-}}$ presents
distinct orientations on $\Sigma_1$.
\end{description}

In face of these previous observations, Theorem \ref{teorema bif
fold-fold} follows straightforward from section \ref{secao fold
fold}.

Note that, as we give the topological behavior of the cases
$\lambda<0$, $\lambda=0$ and $\lambda>0$ it is easy to construct the
bifurcation diagram of \ref{eq pert sing 2} when $\lambda \in \R$.\\

\medskip

\noindent {\textbf{Acknowledgments.}} The first  author is partially
supported by a FA\-PESP$-$BRA\-ZIL grant 2007/08707-5. This work is
partially realized at UFG/Brazil as a part of project number 35799.


\begin{thebibliography}{99}

%

\bibitem{BST} {\sc C.A. Buzzi, P.R. da Silva and M.A. Teixeira},
{\it Singular perturbation problems for time reversible systems},
Proc. Amer. Math. Soc., \textbf{133} (2005), 3323-3331.

\bibitem{Claudio-PR-Marco} {\sc C.A. Buzzi, P.R. da Silva  and M. A. Teixeira},
{\it A singular approach to discontinuous vector fields on the
plane}, Journal of Differential Equations, \textbf{231} (2006),
633-655.
%

\bibitem{DR} {\sc F. Dumortier and R. Roussarie},
{\it Canard cycles and center manifolds}, Memoirs Amer. Mat. Soc.
\textbf{121}, 1996.

\bibitem{F} {\sc N. Fenichel},
{\it Geometric singular perturbation theory for ordinary
differential equations}, Journal of Differential Equations
\textbf{31} (1979), 53--98.

\bibitem{Fi} {\sc A.F. Filippov},
{\it Differential equations with discontinuous righthand sides},
Mathematics and its Applications (Soviet Series), Kluwer Academic
Publishers-Dordrecht, 1988.




\bibitem{K} {\sc V. S. Kozlova},
{\it Roughness of a discontinuous system}, Vestinik Moskovskogo
Universiteta, Matematika \textbf{5} (1984), 16--20.



\bibitem{LST-Regularization-2006} {\sc J. Llibre, P.R. Silva and M.A. Teixeira},
{\it  Regularization of discontinuous vector fields via singular
perturbation}, J. Dynam. Differential Equation \textbf{19} (2006),
309–-331.



\bibitem{LST} {\sc J. Llibre, P.R. Silva and M.A. Teixeira},
{\it  Sliding vector fields via slow$-$fast systems}, Bulletin of
the Belgian Mathematical Society Simon Stevin \textbf{15-5} (2008),
851--869.




\bibitem{LT-Regularization-1997} {\sc J. Llibre and M.A. Teixeira},
{\it  Regularization of discontinuous vector fields in dimension
three}, Discrete Contin. Dynam. Systems \textbf{3} (1997), 235-–241.






\bibitem{ST} {\sc J. Sotomayor and M.A. Teixeira},
{\it Regularization of  discontinuous vector fields}, International
Conference on Differential Equations, Lisboa (1996), 207--223.



\bibitem{S} {\sc P. Szmolyan},
{\it Transversal heteroclinic and homoclinic orbits in singular
perturbation problems}, Journal of Differential Equations
\textbf{92} (1991), 252--281.


\bibitem{T} {\sc M.A. Teixeira},
{\it  Generic singularities of discontinuous vector fields}, An. Ac.
Bras. Cienc. \textbf{53}, n$^{o}$2, (1991), 257--260.


\end{thebibliography}
\end{document}